\documentclass[10pt]{article}

\usepackage{amsfonts} 
\usepackage{amsmath,hhline,latexsym,mathrsfs}
\usepackage{amssymb}
\usepackage{relsize}
\DeclareMathAlphabet\mathbfcal{OMS}{cmsy}{b}{n}

\usepackage[latin1]{inputenc}

\usepackage{hyperref}

\usepackage{color}
\usepackage{graphicx}

\newtheorem{theorem}{\sc Theorem}[section]
\newtheorem{lemma}{\sc Lemma}[section]
\newtheorem{definition}{\sc Definition}[section]
\newtheorem{proposition}{\sc Proposition}[section]
\newtheorem{corollary}{\sc Corollary}[section]
\newtheorem{remark}{Remark}

\newcommand{\mat}[1]{\mbox{\boldmath{$#1$}}}

\newcommand{\dis}{\displaystyle}
\newcommand{\alphasm}{\mathsmaller{\alpha}}
\newcommand{\eps}{\varepsilon}
\newcommand{\om}{\omega}
\newcommand{\Om}{\Omega}
\newcommand{\ph}{\varphi}

\newcommand{\lambar}{\overline{\lambda}}
\newcommand{\ovec}{\mathbf{0}}
\newcommand{\Avec}{\mathbf{A}}
\newcommand{\Pvec}{\mathbf{P}}
\newcommand{\Bvec}{\mathbf{B}}
\newcommand{\Dvec}{\mathbf{D}}
\newcommand{\Mvec}{\mathbf{M}}
\newcommand{\Fvec}{\mathbf{F}}
\newcommand{\Yvec}{\mathbf{Y}}
\newcommand{\Evec}{\mathbf{E}}
\newcommand{\fvec}{\mathbf{f}}

\newcommand{\mubar}{\overline{\boldsymbol{\mu}}}

\newcommand{\ybar}{\overline{y}}

\newcommand{\Fin}{\hfill$\Box$}

\newcommand{\N}{\mbox{$I \kern -4pt N$}}
\newcommand{\Q}{\mbox{$Q \kern -8pt I$}}
\newcommand{\R}{\mbox{$I \kern -4pt R$}}
\newcommand{\C}{\mbox{$C \kern -8pt I$}}

\newcommand{\jnt}{\dis\int}
\newcommand{\jjntQT}{\jnt\!\!\!\!\jnt_{Q_{T}}}
\newcommand{\jjntqT}{\jnt\!\!\!\!\jnt_{q_{T}}}

\newcommand{\pvec}{\mathbf{p}}

\newcommand{\uvec}{\mathbf{u}}
\newcommand{\Lvec}{\mathbf{L}}
\newcommand{\Gvec}{\mathbf{G}}
\newcommand{\Cvec}{\mathbf{C}}

\def\R{{\bf R}}
\title{Inverse problems for linear parabolic equations using mixed formulations - Part 1 : Theoretical analysis %\thanks{Grants or other notes
}
\author{ 
	\textsc{Arnaud M\"unch}\thanks{Laboratoire de Math\'ematiques, Universit\'e Blaise Pascal (Clermont-Ferrand 2), 
	UMR CNRS 6620, Campus de C\'ezeaux, 63177, Aubi\`ere, France. E-mails: {\tt arnaud.munch@math.univ-bpclermont.fr.}}	\quad 
	\and
	\textsc{Diego A. Souza}\thanks{Dpto.\ EDAN, University of Sevilla, 41080~Sevilla, Spain and
	Departamento de Matem\'{a}tica, Universidade Federal da Para\'iba, 58051-900, Jo\~{a}o Pessoa--PB, Brazil.
	E-mail: {\tt desouza@us.es}. Partially supported by CAPES (Brazil) and grant~MTM2010-15592 (DGI-MICINN, Spain).}}
	
\begin{document}

\maketitle

\begin{abstract}
We introduce in this document a direct method allowing to solve numerically inverse type problems for linear parabolic equations. We consider the reconstruction of the full solution of the parabolic equation posed in $\Omega\times (0,T)$ - $\Omega$ a bounded subset of $\mathbb{R}^N$ - from a partial distributed observation. We employ a least-squares technique and minimize the $L^2$-norm of the distance from the observation to any solution. Taking the parabolic equation as the main constraint of the problem, the optimality conditions are reduced to a mixed formulation involving both the state to reconstruct and a Lagrange multiplier. The well-posedness of this mixed formulation - in particular the inf-sup property - is a consequence of classical energy estimates. We then reproduce the arguments to a linear first order system, involving the normal flux, equivalent to the linear parabolic equation. The method, valid in any dimension spatial dimension $N$, may also be employed to reconstruct solution for boundary observations. 

With respect to the hyperbolic situation considered in \cite{NC-AM-InverseProblems} by the first author, the parabolic situation requires - due to regularization properties -  the introduction of appropriate weights function so as to make the problem numerically stable. 
 
% This paper deals with the numerical computation of null controls for the wave equation with a potential. The goal is to compute approximations of controls that drive the solution from a prescribed initial state to zero at a large enough controllability time. In [\textit{C\^{\i}ndea, Fern\'andez-Cara \& M\"unch, Numerical controllability of the wave equation through primal methods and Carleman estimates, 2013}], a so called primal method is described leading to a strongly convergent approximation of boundary controls : the controls minimize quadratic weighted functionals involving both the control and the state and are obtained by solving the corresponding optimality condition. In this work, we adapt the method to approximate the control of minimal square-integrable norm. The optimality conditions of the problem are reformulated as a mixed formulation involving both the state and its adjoint. We prove the well-posedeness of the mixed formulation (in particular the inf-sup condition) then discuss several numerical experiments. The approach covers both the boundary and the inner controllability. For simplicity, we present the approach in the one dimensional case.

\end{abstract}

%\tableofcontents

\section{Introduction - Inverse problems for linear parabolic equations}
\label{sec:intro}

Let $\Om\subset \mathbb{R}^N$ be a bounded domain set whose boundary $\partial\Om$ is regular enough 
	(for instance of class $C^2$). For any $T>0$, we note $Q_T:=\Omega\times (0,T)$ and $\Sigma_T:=\partial\Omega\times (0,T)$.
	
We are concerned with inverse type problems for the following linear parabolic type equation%
\begin{equation}
\label{eq:heat}
	\left\{
		\begin{array}{lll}
   			y_{t} - \nabla\cdot(c(x) \nabla y) + d(x, t) y= f 				& \textrm{in}& Q_T,    \\
   			y  = 0										  		&\textrm{on}& \Sigma_T, \\
   			y(x, 0) = y_0(x) 								  		& \textrm{in}& \Om.
   		\end{array} 
 	\right.
\end{equation}
%	{\color{green} Arnaud, the diffusion matrix $c$ can depend of $t$. You just must suppose that $c:=(c_{i,j})\in C^0([0,T];C^1(\overline\Om;\mathcal{M}_{N}(\mathbb{R})))$ with $(c(x,t)\xi,\xi)\geq c_0|\xi|^2$ 
%	in~$\overline Q_T~(c_0>0)$}\\
	We assume that $c:=(c_{i,j})\in C^1(\overline\Om;\mathcal{M}_{N}(\mathbb{R}))$ with $(c(x)\xi,\xi)\geq c_0|\xi|^2$ 
	for any $x\in\overline\Om$ and $\xi\in \mathbb{R}^N$~$(c_0>0)$, $d \in L^\infty(Q_T)$ and $y_0 \in L^2(\Om)$; $f = f(x,t)$ is a \textit{source} term (a function 
	in~$L^2(Q_T)$) and $y=y(x,t)$ is the associated state. 

In the sequel, we shall use the following notation :
\begin{equation}
%\label{eq:L}
	L\, y:=y_{t} - \nabla\cdot(c(x) \nabla y) + d(x, t) y, \qquad L^{\star}\ph:=-\ph_{t}- \nabla\cdot(c(x) \nabla \ph)+ d(x,t)\ph. \nonumber
\end{equation}
 
For any $y_0\in L^2(\Omega)$ and $f\in L^2(Q_T)$, there exists exactly one solution $y$ to 
	(\ref{eq:heat}), with the regularity $y \in C^0([0, T]; L^2(\Om)) \cap L^2(0, T; H^{1}_0(\Om))$ (see \cite{HarauxCazenave,LionsMagenes1972}).

Let now $\omega$ be any non empty open subset of $\Omega$ and let $q_T:=\omega \times (0,T)\subset Q_T$.
A typical inverse problem for (\ref{eq:heat}) (see \cite{isakov}) is the following one~:  from an \textit{observation} or \textit{measurement} $y_{obs}$ in $L^2(q_T)$ on the open set $q_T$, we want to recover a solution $y$ of the boundary value problem (\ref{eq:heat}) which coincides with the observation on $q_T$. Introducing the operator $P: \mathcal{Y} \to L^2(Q_T)\times L^2(q_T)$ defined by $Py:=(Ly,y_{\vert q_T})$ where the space $\mathcal{Y}$ is an appropriate Hilbert space (defined in Section \ref{sec1_direct}), the problem is reformulated as\,: 
\begin{equation}
\label{IP}\tag{$IP$}
	\text{\it find } y\in \mathcal{Y} \text{ \it solution of } Py=(f,y_{obs}). 
\end{equation}

From the unique continuation property for (\ref{eq:heat}), if $y_{obs}$ is a restriction to $q_T$ of a solution of (\ref{eq:heat}), then the problem is well-posed in the sense that the state $y$ corresponding to the pair $(y_{obs},f)$ is unique, i.e. $P$ is a 
bijective operator from $\mathcal{Y}$ to its range $\mathcal{R}(P)$.

%On the other hand, if the source $f$ is not fixed, then in general, the solution $y$ is not unique. Consider for instance a source $f$ supported in a set near $\Omega\times \{T\}$ and disjoint from $q_T$; from the finite speed of propagation, the source $f$ will not affect the solution $y$ in $q_T$. The uniqueness is recover adding some specific property on $f$. 

%We address in this work the numerical approximation of the problem $(IP)$ with the aim to produce strongly convergent family of approximation. 

In view of the unavoidable uncertainties on the data $y_{obs}$ (coming from measurements, numerical approximations, etc), Problem (\ref{IP}) needs to be relaxed. In this respect, the most natural (and widely used in practice) approach consists in introducing the following extremal problem (of least-squares type)
\begin{equation}
\label{extremal_problem} \tag{\it LS}
\quad
\left\{
\begin{aligned}
&\textrm{minimize over } \mathcal{H} \quad J(y_0):=\frac{1}{2}  \Vert \rho_0^{-1}(y-y_{obs}) \Vert^2_{L^2(q_T)}\\
& \textrm{where} \quad y \quad \textrm{solves} \quad (\ref{eq:heat}),
\end{aligned}
\right.
\end{equation}
since $y$ is uniquely and fully determined from the data $y_0$. $\rho_0$ denotes a appropriate positive weight while $\mathcal{H}$ denotes a Hilbert space related to the space $\mathcal{Y}$: roughly, $\mathcal{H}$ is the set of all initial data $y_0$ for which the solution of (\ref{eq:heat}) satisfies $\rho_0^{-1}y\in L^2(q_T)$. 

Here the constraint $y-y_{obs}=0$ in $L^2(q_T)$ is relaxed; however, if $y_{obs}$ is a restriction to $q_T$ of a solution of (\ref{eq:heat}), then problems \eqref{extremal_problem} and \eqref{IP} coincide. A minimizing sequence for $J$ in $\mathcal{H}$ is easily defined in term of the solution of an auxiliary adjoint problem. However, on a numerical point of view, this extremal problem has mainly two independents drawbacks:
\begin{itemize}
\item First, it is in general not possible to minimize over a discrete subspace of the set $\{y; Ly-f=0\}$ subject to the equality (in $L^2(Q_T)$) $Ly-f=0$. Therefore, the minimization procedure first requires the discretization of the functional $J$ and of the system (\ref{eq:heat}); this raises the issue, when one wants to prove some convergence result of any discrete approximation, of the uniform coercivity property (typically here some uniform discrete observability inequality for the adjoint solution) of the discrete functional with respect to the approximation parameter. As far as we know, this delicate issue has received answers only for specific and somehow academic situations (uniform Cartesian approximation of $\Omega$, constant coefficients in \eqref{eq:heat}, etc). We refer to \cite{boyer_canum2012,munch-zuazua}.
\item Second, in view of the regularization property of the heat kernel, the space of initial data $\mathcal{H}$ for which the corresponding solution of (\ref{eq:heat}) belong to $L^2(q_T)$ is a huge space. Its contains in particular the negative Sobolev space $H^{-s}(\Omega)$ for any $s>0$ and therefore is very hard to approximate numerically. For this reason, the reconstruction of the initial condition $y_0$ of (\ref{eq:heat}) from a partial observation in $L^2(q_T)$ is therefore known to be numerically severally ill-posed and requires, within this framework, a regularization to enforce that the minimizer belongs, for instance, to $L^2(\Omega)$ much easier to approximate (see \cite{engl}). The situation is analogous for the so-called backward heat problem, where the observation on $q_T$ is replaced by a final time observation. We refer to (\cite{choulli-yamamoto,puel_rio,wang-yamamoto}) where this ill-posedness is discussed.    
\end{itemize}

The main reason of this work is to reformulate problem (\ref{extremal_problem}) and show that the use of variational methods may overcome these two drawbacks. 

Preliminary, we also mention that the quasi-reversibility method initially introduced in \cite{lattes_lions69} may be employed to address problem (\ref{IP}). This kind of methods, which falls into the category of regularization methods, reads as follows:  for any $\eps>0$, find $y_{\eps}\in \mathcal{Y}$ the solution of 
\begin{equation}\label{FV}\tag{$QR$}
 \langle Py_\eps,P\overline{y} \rangle_{L^2(Q_T)\times L^2(q_T)} + \eps \langle y_{\eps},\overline{y} \rangle_{\mathcal{Y}} = \left\langle(f,y_{obs}), P\overline{y}\right\rangle_{L^2(Q_T)\times L^2(q_T),L^2(Q_T)\times L^2(q_T)},
\end{equation}
for all $\overline{y}\in \mathcal{Y}$, 
where $\eps>0$ is a Tikhonov like parameter which ensures the well-posedness. We refer to the book \cite{Klibanov-book} (and the references therein) and more recently to \cite{bourgeois_parabolic} where the Cauchy problem for the heat equation is addressed. Remark that (\ref{FV}) can be viewed as a least-squares problem since the solution $y_{\eps}$ minimizes over $\mathcal{Y}$ the functional $y\to \Vert P y - (f,y_{obs})\Vert^2_{L^2(Q_T)\times L^2(q_T)}+ \eps \Vert y\Vert^2_{\mathcal{Y}}$. Eventually, if $y_{obs}$ is a restriction to $q_T$ of a solution of (\ref{eq:heat}), 
the corresponding $y_{\eps}$ converges in $L^2(Q_T)$ toward to the solution of \eqref{IP} as $\eps\to 0$. There, unlike in Problem \eqref{extremal_problem}, the unknown is the state variable $y$ itself (as it is natural for elliptic equations) so that any standard numerical methods based on a conformal approximation of the space $\mathcal{Y}$ together with appropriate global observability inequalities allow to obtain a convergent approximation of the solution. In particular, there is no need to prove discrete observability inequalities. We refer to the book \cite{BeilinaKlibanov14}. 

On the other hand, we also mention that Luenberger type observers approach, recently used to address reconstruction problems for hyperbolic equation by exploiting reversibility properties (see \cite{cindea_moireau,ramdani2012}) are a priori un-effective in the parabolic situation given by (\ref{eq:heat}).

In the spirit of the works \cite{bourgeois_parabolic, Klibanov-book}, we explore the direct resolution of the optimality conditions associated to the extremal problem \eqref{extremal_problem}, without Tikhonov parameter while keeping the $y$ as the unknown of the problem. This strategy, advocated in \cite{puel_rio}, avoids any iterative process and allows a stable numerical framework: it has been successfully applied in the closely related context of the exact controllability of (\ref{eq:heat}) in \cite{EFC-AM-sema,munch_desouza} and also to inverse problems for hyperbolic equations in \cite{NC-AM-InverseProblems, NC-AM-InverseProblemsBoundary}. Keeping $y$ as the main variable, the idea is to take into account the state constraint $Ly-f=0$ with a Lagrange multiplier. This allows to derive explicitly the optimality systems associated  to \eqref{extremal_problem} in term of an elliptic mixed formulation and therefore reformulate the original problem. Well-posedness of such new formulation is related to classical energy estimates and unique continuation properties while the stability is guarantees by some global observability inequality for the homogeneous parabolic equation. 

The outline of this paper is as follow. In Section \ref{recovering_y}, we consider the least-squares problem \eqref{P} and reconstruct the solution of the parabolic equation from a partial observation localized on a subset $q_T$ of $Q_T$. For that, in Section \ref{sec1_direct}, we associate to \eqref{P} the equivalent mixed formulation (\ref{eq:mf}) which relies on the optimality conditions of the problem. Using the unique continuation for the equation (\ref{eq:heat}), we show the well-posedness of this mixed formulation, in particular, we check the Babuska-Brezzi inf-sup condition (see Theorem \ref{th:mf}). Interestingly, in Section \ref{sec1_dual}, we also derive a equivalent dual extremal problem, which reduces the determination of the state $y$ to the minimization of a elliptic functional with respect to the Lagrange multiplier. Then, in Section \ref{recovering_y_fo}, we adapt these arguments to the first order mixed system (\ref{eq:heat_mixed}), equivalent to the parabolic equation. There, the flux variable $\mathbf{p}:=c(x)\nabla y$ appears explicitly in the formulation and allows to reduce the order of regularity of the involved functional spaces. The underlying inf-sup condition is obtained by adapting a Carleman inequality due to Imanuvilov, Puel and Yamamoto (see \cite{imanuv_puel_yama_2}).
The existence and uniqueness of weak solution to this first order system is studied in the Appendix.
Section \ref{sec_conclusion} concludes with some remarks and perspectives: in particular, we highlight why the mixed formulations developed and analyzed here are suitable at the numerical level to get a robust approximation of the variable $y$ on the whole domain $Q_T$.

%%% NEW SECTION %%% 

\section[Second order mixed formulation]{Recovering the solution from a partial observation: a second order mixed formulation}\label{recovering_y}

In this section, assuming that the initial data $y_0$ is unknown, we address the inverse problem \eqref{IP}. Without loss of generality, in view of the linearity of (\ref{eq:heat}), we assume that the source term $f$ is zero: $f\equiv 0$ in $Q_T$.

\subsection{Direct approach : Minimal local weighted $L^2$-norm: a first mixed formulation }
\label{sec1_direct}

Let $\rho_{\star}\in \mathbb{R}^+_{\star}$ and let $\rho_0\in\mathcal{R}$ with
\begin{equation}
	\mathcal{R}:=\{w: w\in C(Q_T); w\geq \rho_{\star}>0 \,\, \textrm{in}\,\, Q_T; w\in L^{\infty}(\Om\times(\delta,T)) \,\, \forall \delta>0\}   \label{def_spaceR}
\end{equation}
so that in particular, the weight $\rho_0$ may blow up as $t\to 0^{+}$. 

We define the space 
\[
	 \mathcal{Y}_0 := \left\{ y \in C^2(\overline Q_T): y = 0 \textrm{ on }  \Sigma_T\right\}
\]
and for any $\eta>0$ and any $\rho\in \mathcal{R}$, the bilinear form by 
\[
	(y, \ybar)_{\mathcal{Y}_0} := \iint_{q_T} \rho^{-2}_0 y \,\ybar\, dx\,dt  + \eta\iint_{Q_T} \rho^{-2}L\, y L\, \ybar\, dx\,dt, 
	\quad \forall y, \ybar \in \mathcal{Y}_0.
\]
The introduction of the weight $\rho$ which does not appear in the original problem will be motivated at the end of this section.
From the unique continuation property for (\ref{eq:heat}), this bilinear form defines for any $\eta>0$ a scalar product. 

Let then $\mathcal{Y}$ be the completion of the space $\mathcal{Y}_0$ for this scalar product. We denote the norm over 
$\mathcal{Y}$ by $\Vert\cdot\Vert_\mathcal{Y}$ such that
\begin{equation}\nonumber
\Vert y\Vert^2_\mathcal{Y}:= \Vert \rho_0^{-1} y\Vert^2_{L^2(q_T)}  + \eta \Vert \rho^{-1}L\,y\Vert^2_{L^2(Q_T)}, \quad \forall y\in\mathcal{Y}.
\end{equation}

Finally, we define the closed subset $\mathcal{W}$ of $\mathcal{Y}$ by 

\begin{equation}\nonumber
\mathcal{W}:=\biggl\{ y\in \mathcal{Y} : \rho^{-1}Ly=0\; \textrm{in}\; L^2(Q_T)   \biggr\}
\end{equation}
and we endow $\mathcal{W}$ with the same norm than $\mathcal{Y}$.

We then define the following extremal problem : 

\begin{equation}
\label{P}\tag{$P$}
	\left\{
		\begin{array}{l}
			\dis \hbox{Minimize }\ J(y) := {1 \over 2} \iint_{q_T} \rho^{-2}_0 |y(x,t)-y_{obs}(x,t)|^2 \,dx\,dt \\
			\noalign{\smallskip}
			\hbox{Subject to }\ y \in \mathcal{W}.
		\end{array}
	\right.
\end{equation}

This extremal problem is well-posed: the functional $J$ is continuous, strictly convex and is such that $J(y)\to +\infty$ as $\Vert y\Vert_\mathcal{Y}\to+ \infty$. Note also that the solution of the problem in $\mathcal{W}$ does not depend on $\eta$ nor $\rho$.
Moreover, for any $y\in \mathcal{W}$, $Ly=0$ a.e. in $Q_T$ and $\Vert y\Vert_\mathcal{Y}=\Vert \rho_0^{-1}y\Vert_{L^2(q_T)}$ so that the restriction $y(\cdot,0)$ belongs by definition to the abstract space $\mathcal{H}$: consequently, extremal problems (\ref{extremal_problem}) and (\ref{P}) are equivalent. 

In order to solve problem (\ref{P}), we have to deal with the constraint equality $\rho^{-1}Ly=0$ which appears in  
$\mathcal{W}$. Proceeding as in \cite{NC-AM-InverseProblems,munch_desouza}, we introduce a Lagrange multiplier and the following mixed formulation: find $(y,\lambda)\in \mathcal{Y}\times L^2(Q_T)$ solution of 
\begin{equation} 
\label{eq:mf}
	\left\{
		\begin{array}{rcll}
			\noalign{\smallskip} 
			a(y, \ybar) +  b(\ybar,\lambda) & = & l(\ybar) & \quad \forall \ybar \in   \mathcal{Y},\\
			\noalign{\smallskip}
			                          b( y, \lambar) & = & 0 & \quad \forall \lambar \in L^2(Q_T),
		\end{array}
	\right.
\end{equation}
	where
\begin{align}
	\nonumber
	& a : \mathcal{Y} \times\mathcal{Y} \to \mathbb{R},  \quad  a(y, \ybar)  := 	  \jjntqT\rho_0^{-2}y \,\ybar\, dx\,dt\\
	\nonumber
	& b: \mathcal{Y} \times L^2(Q_T)\to \mathbb{R}, \quad  b(y, \lambda) := \jjntQT\rho^{-1}Ly \,\lambda\,dx\,dt\\
	\nonumber
	& l:\mathcal{Y}\to \mathbb{R},  \quad  l(y) :=\jjntqT\rho^{-2}_0 y \,y_{obs}\, dx\,dt .
\end{align}

	We have the following result : 

\begin{theorem}\label{th:mf}
	Let $\rho_0\in \mathcal{R}$ and $\rho\in \mathcal{R}\cap L^{\infty}(Q_T)$. 
\begin{enumerate}
	\item 
		The mixed formulation \eqref{eq:mf} is well-posed.
	\item 
		The unique solution $(y,\lambda)\in \mathcal{Y}\times L^2(Q_T)$ is the unique saddle-point of the Lagrangian 
		 ${\mathcal{L} :  \mathcal{Y}\times L^2(Q_T) \to \mathbb{R}}$ defined by
\begin{equation}
\label{eq:calLtilde}
	%\mathcal{ L}(y, \lambda) = {1 \over 2} \iint_{q_T} \rho^{-2}_0|y(x,t)-y_{obs}(x,t)|^2 \,dx\,dt+\iint_{Q_T}\rho^{-1}Ly \lambda\,dx\,dt.
	\mathcal{L}(y,\lambda):=\frac{1}{2}a(y,y)+b(y,\lambda)-l(y).
\end{equation}
	\item 
		The solution $(y,\lambda)$ satisfies the estimates
\begin{equation}
\label{estimative_sol_weig}
	\|y\|_\mathcal{Y}\leq\|\rho_0^{-1}y_{obs}\|_{L^2(q_T)},
	\quad \|\lambda\|_{L^2(Q_T)}\leq2\sqrt{\rho_{\star}^{-2}\Vert \rho\Vert^2_{L^{\infty}(Q_T)}+ \eta}\,\|\rho_0^{-1}y_{obs}\|_{L^2(q_T)}.
\end{equation}
\end{enumerate}
\end{theorem}
\par\noindent
\textsc{Proof-} We use classical results for saddle point problems  (see \cite{brezzi_new}, chapter 4).

We easily obtain the continuity of the symmetric and positive bilinear form $ a$ over 
	$\mathcal{Y}\times \mathcal{Y}$, the continuity of the bilinear form $b$ over $\mathcal{Y}\times L^2(Q_T)$ and the continuity of the linear form $l$ over 
	$\mathcal{Y}$.  In particular, we get
\begin{equation}\label{norm_linear_form_weig}
	\|l\|_{ \mathcal{Y}'}\leq\|\rho^{-1}_0y_{obs}\|_{L^2(q_T)},\quad \|a\|_{\mathscr{L}^2(\mathcal{Y})}\leq 1,
	\quad  \|b\|_{\mathscr{L}^2(\mathcal{Y},L^2(Q_T))}\leq\eta^{-1/2},
\end{equation}
	where $\mathscr{L}^2(E,F)$ denotes the space of the continuous bilinear functions defined on the product Banach spaces $E\times F$; when 
	$E=F$, we simply write $\mathscr{L}^2(E)$.

Moreover, the kernel $\mathcal{N}( b):=\{y\in \mathcal{Y} : b(y,\lambda)=0\,~ \forall \lambda\in L^2(Q_T)\}$ coincides with
$\mathcal{W}$: we have
$$
a(y,y)=\Vert y\Vert^2_\mathcal{Y}, \qquad \forall y\in \mathcal{N}(b):=\mathcal{W}
$$ 
leading to the coercivity of $a$ over the kernel of $b$. 

Therefore, in view of \cite[Theorem 4.2.2]{brezzi_new}, it remains to check the so-called inf-sup property: there exists $\delta>0$ such that
\begin{equation}
\label{eq:infsup3}
	\inf_{\lambda \in L^2(Q_T)} \sup_{y\in \mathcal{Y}} \frac{ b(y, \lambda)}{\|y\|_\mathcal{Y} \|\lambda\|_{L^2(Q_T)}} \geq \delta.
\end{equation}

We proceed as follows. For any fixed $\lambda^0\in L^2(Q_T)$, using the fact that $\rho$ is bounded in $Q_T$, we define the unique element $y^0$ solution of 
\begin{equation}
	\rho^{-1}L \,y^0=\lambda^0 \,\textrm{ in }\,Q_T, \quad y^0=0\,\textrm{ on }\,\Sigma_T, \quad y^0(\cdot,0)=0 \,\, \textrm{in}\,\, \Omega \nonumber.
\end{equation}

Using energy estimates, we have 
\begin{equation}
\Vert \rho_0^{-1}y^0\Vert_{L^2(q_T)} \leq \rho_{\star}^{-1} \Vert y^0\Vert_{L^2(Q_T)}  \leq \rho_{\star}^{-1}\Vert \rho \lambda^0\Vert_{L^2(Q_T)} 
\leq \rho_{\star}^{-1}\Vert \rho\Vert_{L^{\infty}(Q_T)} \Vert \lambda^0\Vert_{L^2(Q_T)}   \label{estimate_L2L2}
\end{equation}
	which proves that $y^0\in \mathcal{Y}$ and that  
\[
	\sup_{y \in \mathcal{Y}} \frac{b(y, \lambda^0)}{\|y\|_\mathcal{Y} \|\lambda^0\|_{L^2(Q_T)}} 
	\geq \frac{ b(y^0, \lambda^0)}{\|y^0\|_\mathcal{Y} \|\lambda^0\|_{L^2(Q_T)}} 
	= \frac{\|\lambda^0\|_{L^2(Q_T)}}{\left( \|\rho^{-1}_0y^0\|^2_{L^2(q_T)} + \eta \|\lambda_0\|^2_{L^2(Q_T)} \right)^\frac{1}{2}}.
\]

	Combining the above two inequalities, we obtain
\[
	\sup_{y\in \mathcal{Y}} \frac{b(y, \lambda_0)}{\|y\|_\mathcal{Y} \|\lambda_0\|_{L^2(Q_T)}} 
	\geq \frac{1}{\sqrt{\rho_{\star}^{-2}\Vert \rho\Vert^2_{L^{\infty}(Q_T)}  + \eta}}
\]
	and, hence, \eqref{eq:infsup3} holds with $\delta = \left( \rho_{\star}^{-2}\Vert \rho\Vert^2_{L^{\infty}(Q_T)}+ \eta \right)^{-1/2}$.

	The point $(ii)$ is due to the positivity and symmetry of the form $a$.

	The point $(iii)$ is a consequence of classical estimates (see \cite{brezzi_new}, Theorem 4.2.3)\,:
\[
	\|y\|_\mathcal{Y}\leq{1\over\alpha_0}\|l\|_{\mathcal{Y}'},
	\quad\|\lambda\|_{L^2(Q_T)}\leq{1\over\delta}\left(1+{\|a\|_{\mathscr{L}^2(\mathcal{Y})}\over\alpha_0}\right)\|l\|_{\mathcal{Y}^{\prime}},
\]
	where
\[
	\alpha_0:=\inf\limits_{y\in \mathcal{N}(b)} {a(y,y)\over\|y\|_\mathcal{Y}^2}.
\]
Estimates \eqref{norm_linear_form_weig} and the equality $\alpha_0=1$ lead to the results.
	
\Fin

	In order to get a global estimate of the reconstructed solution, we now recall the following important result.
\begin{proposition}[Lemma $3.1$ 	in \cite{cara_guerrero}]\label{carleman_estimate}
Let the weights $\rho_c,\rho_{c,0} \in \mathcal{R}$ $($see \eqref{def_spaceR}$)$ be defined  as follows\,:
 \begin{equation}
\label{weights}
\begin{aligned}
& \rho_c(x,t):=  \exp\left({\beta(x) \over t}\right), \quad \beta(x):=  K_{1}\left(e^{K_{2}} \!-\! e^{\beta_{0}(x)}\right), \\
& \rho_{c,0}(x,t):=  t^{3/2}\rho_c(x,t), \quad \rho_{c,1}(x,t):=  t^{1/2}\rho_c(x,t) 
\end{aligned}
\end{equation}
with $\beta_0 \in C^\infty(\overline{\Omega})$ and where the positive constants $K_{i}$ are sufficiently large (depending on $T$, $c_0$, $\|c\|_{C^1(\overline{\Omega})}$ and $\Vert d\Vert_{L^{\infty}(Q_T)}$) such that 
\[
\beta > 0\quad\textrm{in}\quad \Omega,\, \beta=0\quad\textrm{on}\quad \partial\Omega, \quad \nabla\beta(x)\neq0\quad\forall x\in\overline{\Omega}\setminus\omega.
\]
Then, there exists a constant $C>0$, depending only on $\omega, T$,  such that 
\begin{equation}
\Vert \rho_{c,0}^{-1}y\Vert_{L^2(Q_T)}+\Vert \rho_{c,1}^{-1}\nabla y\Vert_{L^2(Q_T)} \leq C \Vert y\Vert_{\mathcal{Y}_c} \quad \forall y\in \mathcal{Y}_c,   \label{crucial_estimate_est}
\end{equation}
	where $\mathcal{Y}_c$ is the completion of $\mathcal{Y}_{0,c}:=\mathcal{Y}_0$  with respect to the scalar product
\[
	(y, \ybar)_{\mathcal{Y}_{0,c}} = \jjntqT \rho^{-2}_{c,0} y \,\ybar\, dx\,dt  + \eta\jjntQT \rho^{-2}_cLy\, L\ybar\, dx\,dt.
\]
\end{proposition}

The estimate \eqref{crucial_estimate_est} is a consequence of the celebrated global Carleman inequality satisfied by the solution of 
\eqref{eq:heat}, introduced and popularized in \cite{FursikovImanuvilov}. This result implies the following stability estimate which allows to estimate a global norm of the solution $y$ in term of the norm $\mathcal{Y}$.

\begin{corollary}
Let $\rho_0\in \mathcal{R}$ and $\rho\in \mathcal{R}\cap L^{\infty}(Q_T)$ and assume that there exists a positive constant $K$ such that 
\begin{equation}
\rho_0\leq K \rho_{c,0}, \quad \rho\leq K \rho_c \quad \textrm{in}\quad Q_T.   \label{hypK}
\end{equation}
If $(y,\lambda)$ is the solution of the mixed formulation (\ref{eq:mf}), then there exists $C>0$ such that 
\begin{equation}
\Vert\rho_{c,0}^{-1}y\Vert_{L^2(Q_T)}  \leq  C\Vert y\Vert_\mathcal{Y}.  
 \label{crucial_estimate_good}
\end{equation}
\end{corollary}
\par\noindent
\textsc{Proof-} The hypothesis (\ref{hypK}) implies that $\mathcal{Y}\subset \mathcal{Y}_c$. Therefore, estimate (\ref{crucial_estimate_est}) implies that 

\begin{equation}\nonumber
\Vert \rho_{c,0}^{-1}y\Vert_{L^2(Q_T)} \leq C \Vert y\Vert_{\mathcal{Y}_c}\leq C \Vert y\Vert_\mathcal{Y}\leq C 
\Vert\rho_0^{-1} y_{obs}\Vert_{L^2(q_T)}.  
% \label{crucial_estimate_good}
\end{equation}

\begin{remark}
The well-posedness of the mixed formulation (\ref{eq:mf}), precisely the inf-sup property (\ref{eq:infsup3}), is open in the case where the weight $\rho$ is simply in $\mathcal{R}$: in that case, the weight may blow up at time $t=0$. In order to get (\ref{eq:infsup3}), it suffices to prove that the function $z:=\rho_0^{-1}y$ solution of the boundary value problem 
\[
\rho^{-1}L(\rho_0 z)=\lambda^0 \textrm{ in } Q_T, \quad z=0 \textrm{ on } \Sigma_T, \quad z(\cdot,0)=0\textrm{ in }\Omega
\]
for any $\lambda^0\in L^2(Q_T)$ satisfies the following estimate for some positive constant $C$
\[
\Vert z\Vert_{L^2(q_T)} \leq C \Vert \rho^{-1}L(\rho^0 z)\Vert_{L^2(Q_T)}.
\]
In the cases of interest for which both $\rho_0$ and $\rho$ blow up at $t\to 0^+$ (for instance given that $\rho_{c,0}$ and $\rho_c$), this estimate is open and does not seem to be a consequence of the estimate 
(\ref{crucial_estimate_est}).
\end{remark}

Let us now comment the introduction of the weight $\rho_0$ in the problem (\ref{P}). The space $\mathcal{Y}_c$, which contains the element $y$ such that $\rho_c^{-1}Ly\in L^2(Q_T)$ and $\rho_{c,0}^{-1}y\in L^2(q_T)$ satisfies the embedding $\mathcal{Y}_c\subset C^0([\delta,T],H_0^1(\Omega))$ for any $\delta>0$ (see \cite{cara_guerrero}). Under the condition (\ref{hypK}), the same embedding holds for $\mathcal{Y}$. In particular, there is no control of the restriction of the solution at time $t=0$, which is of course due to the regularization property of the heat kernel. Consequently, from the observation $y_{obs}\in L^2(q_T)$ and the knowledge of $Ly\in L^2(Q_T)$, there is no hope to recover - for a Sobolev norm - the solution of $y$ at the initial time $t=0$. It is then suitable to add to the cost $J$, a weight $\rho_0^{-1}$ 
that vanishes at time $0$.  The weight $\rho$ is introduced here for similar reasons. Remark that the solution $y$ of (\ref{eq:mf}) belongs to $\mathcal{W}$ and therefore does not depend on $\rho$ (recall that $\rho$ is strictly positive): this is in agreement with the fact that $\rho$ does not appear in the equivalent problem (\ref{extremal_problem}). However, very likely, a singular behavior for the $L^2(Q_T)$ function $Ly$ occurs as well near $\Omega\times \{0\}$ so that the constraint $Ly=0$ in $L^2(Q_T)$ is too "strong" and must be replaced - for numerical purposes -  by the relaxed one $\rho^{-1}Ly=0$ in $L^2(Q_T)$ with $\rho^{-1}$ small near $\Omega\times \{0\}$. Remark that this is actually the effect and the role of the Carleman type weight $\rho^c$ defined in (\ref{weights}). As a partial conclusion, the introduction of appropriate weights in the cost $J$ allows to use the estimate (\ref{crucial_estimate_good}) and to guarantee a Lipschitz stable reconstruction of the solution $y$ on the whole domain except at the initial time. 

We also emphasize that the mixed formulation (\ref{eq:mf}) is still well defined with constant weights, $\rho$ and $\rho_0$ equals to one, but leading to weaker stability estimates and reconstruction results. We refer to \cite{choulli_book,choulli-yamamoto}.

\

Furthermore, - at  the numerical level - it is also very convenient to ``augment" the Lagrangian  $\mathcal{L}$ (see \cite{fortinglowinski}) and consider instead the Lagrangian $\mathcal{L}_r$ defined for any $r\geq 0$ by 
\begin{equation}
\left\{
\begin{aligned}
& \mathcal{L}_r(y,\lambda):=\frac{1}{2}a_r(y,y)+b(y,\lambda)-l(y), \\
& a_r(y,y):=a(y,y)+r\Vert \rho^{-1} Ly\Vert^2_{L^2(Q_T)}. 
\end{aligned}
\right.
 \nonumber
\end{equation}
Since $a_r(y,y)=a(y,y)$ on $\mathcal{W}$, the Lagrangian $\mathcal{L}$ and $\mathcal{L}_r$ share the same saddle-point. The nonnegative number $r$ is an augmentation parameter. 

\begin{remark}\label{weak_lambda}
	The first equation of the mixed formulation \eqref{eq:mf} reads as follows: 
\begin{equation}
\iint_{q_T}\rho_0^{-2}\,y \,\overline y\, dx\,dt +\iint_{Q_T} \rho^{-1}L\overline y\, \lambda\, dx\, dt = \iint_{q_T} \rho_0^{-2}\,y_{obs} \,\overline{y} \, dx\,dt \quad \forall \overline y\in \mathcal{Y}.
\end{equation}

	But this means that $\rho^{-1}\lambda\in L^2(Q_T)$ is solution of the parabolic equation in the transposition sense, i.e.
	$\rho^{-1}\lambda$ solves the problem\,: 
\begin{equation}
\label{eq:system_lambda}
	\left\{
		\begin{array}{lcl}
   			L^{\star}\,(\rho^{-1}\lambda)=-\rho_0^{-2}\,(y-y_{obs})1_\om	&\textrm{in}&	Q_T,    	\\
   			\rho^{-1}\lambda  = 0									&\textrm{on}&	\Sigma_T,	\\
			(\rho^{-1}\lambda)(\cdot,T)=0							&\textrm{in}&	\Om.
   		\end{array} 
 	\right.
\end{equation}
where $1_\om$ denotes the characteristic function associated to the open set $\omega$.

	Therefore, $\rho^{-1}\lambda$ (defined in the weak sense) is the solution of a backward parabolic problem with zero initial state and right hand side  $-\rho_0^{-2}(y-y_{obs})1_\om$ in $L^2(q_T)$. This implies in particular that $\rho^{-1}\lambda$ is more regular than $L^2(Q_T)$: precisely, one can see that  $\rho^{-1}\lambda$ belongs to
	$C^0([0,T];H^1_0(\Omega))\cap L^2(0,T;H^2(\Omega)\cap H_0^1(\Omega))$.
	
\begin{itemize}
\item[\textbullet]  Moreover, if $y_{obs}$ is the restriction to $q_T$ of a solution of (\ref{eq:heat}), then the unique multiplier $\lambda$, solution of (\ref{eq:system_lambda}), must vanish almost everywhere. In that case, 
we have ~ $\sup_{\lambda\in L^2(Q_T)}\inf_{y\in \mathcal{Y}} \mathcal{L}_r(y,\lambda) 
= \inf_{y\in \mathcal{Y}} \mathcal{L}_r(y,0)=\inf_{y\in \mathcal{Y}} J_r(y)$ with 
\begin{equation}
\label{def_Jyr}
J_r(y):=\frac{1}{2}\Vert \rho_0^{-1}\,(y-y_{obs})\Vert^2_{L^2(Q_T)} + \frac{r}{2}\Vert \rho^{-1}Ly\Vert^2_{L^2(Q_T)}.
\end{equation}
The corresponding variational formulation, well-posed for $r>0$ is then : find $y\in \mathcal{Y}$ such that 
\begin{equation}
a_r(y,\overline{y})= l(\overline{y}), \quad \forall \overline{y}\in \mathcal{Y}.  \nonumber
\end{equation}
\item[\textbullet] In the general case, the mixed formulation can be rewritten as follows: find $(y,\lambda)\in 
\mathcal{Y}\times L^2(Q_T)$ solution of 
\begin{equation}
\label{QRbis}
\hspace{-0.4cm}
\left\{
\begin{aligned}
\!\langle P_r y, P_r \overline{y}\rangle_{L^2(Q_T)\times L^2(q_T)} + \!\langle \rho^{-1}L\overline{y},\lambda \rangle_{L^2(Q_T),L^2(Q_T)}   & \!=\! \langle (0,\rho^{-1}_0y_{obs}),P_r \overline{y} \rangle_{L^2(Q_T)\times L^2(q_T)}~\, \forall \overline{y}\in \mathcal{Y}, \\
\langle \rho^{-1}Ly, \overline{\lambda} \rangle_{L^2(Q_T),L^2(Q_T)}   & = 0 \quad \forall \overline{\lambda}\in L^2(Q_T)
\end{aligned}
\right.
\end{equation}
with $P_r y:=(\sqrt{r} \rho^{-1}L\,y,\rho_0^{-1}y_{\vert q_T})$. Formulation (\ref{QRbis}) may be seen as generalization of the \eqref{FV} problem (see (\ref{FV}) in the introduction), where the variable $\lambda$ is adjusted automatically (while the choice of the Tikhonov type parameter $\eps$ in (\ref{FV}) is in general a delicate issue).
\end{itemize}	
\end{remark}

\

The optimality system \eqref{eq:system_lambda} can be used to define a equivalent saddle-point formulation, very suitable at the 
	numerical level.  Precisely, we introduce - in view of \eqref{eq:system_lambda} - the space $\Lambda$ by 
\[
	\Lambda:=\{\lambda: \rho^{-1}\lambda\in C^0([0,T]; L^2(\Omega)) ,\,  
	\rho_0\,L^\star(\rho^{-1}\lambda)\in L^2(Q_T),\,\rho^{-1}\lambda=0\hbox{ on }\Sigma_T,\,(\rho^{-1}\lambda)(\cdot,T)=0\}.
\]
	Endowed with the scalar product  $\langle \lambda,\overline{\lambda}\rangle_{\Lambda} := 
	\jjntQT (\rho^{-2}\lambda\,\overline{\lambda}+\rho^2_0L^\star(\rho^{-1}\lambda) L^\star(\rho^{-1}\overline{\lambda}))\, dxdt$,
	we first check that %, in view of assumption \eqref{iobs},
	  $\Lambda$ is a Hilbert space.  Then, for any parameter $\alpha\in (0,1)$, we consider the following mixed formulation\,: 
	  find $(y,\lambda)\in \mathcal{Y}\times \Lambda$ such that 
\begin{equation} \label{eq:mfalpha}
	\left\{
		\begin{array}{rcll}
			\noalign{\smallskip}\dis
			 a_{r,\alpha}(y, \overline{y}) + b_{\alpha}(\overline{y}, \lambda) & = & l_{1,\alpha}(\overline{y}), 
			 & \qquad \forall \overline{y} \in \mathcal{Y}, \\
			 \noalign{\smallskip}\dis
			 b_{\alpha}(y, \overline{\lambda})-c_{\alpha}(\lambda,\overline{\lambda})&=&l_{2,\alpha}(\overline{\lambda}),
			 & \qquad \forall \overline{\lambda} \in \Lambda,
		\end{array}
	\right.
\end{equation}
	where
\begin{align}
	& a_{r,\alpha} : \mathcal{Y} \times \mathcal{Y} \to \mathbb{R},  
	\quad 
	a_{r,\alpha}(y,\overline{y}):=(1-\alpha)\jjntqT\rho_0^{-2}y\overline{y}\,dxdt+r\jjntQT\rho^{-2}LyL\overline{y}\,dxdt,
	\nonumber\\
	& b_{\alpha}: \mathcal{Y}\times \Lambda \to \mathbb{R},  
	\quad 
	b_{\alpha}(y,\lambda) := \jjntQT  \rho^{-1}Ly\lambda dt-\alpha \jjntqT  L^\star(\rho^{-1}\lambda)\,y\,dxdt, 
	\nonumber\\ 
	& c_{\alpha}: \Lambda \times \Lambda  \to \mathbb{R},  
	\quad 
	c_{\alpha}(\lambda, \overline{\lambda}) := \alpha\jjntQT \rho^{2}_0L^\star(\rho^{-1}\lambda)\,
	L^\star(\rho^{-1}\overline{\lambda}) \, dxdt, \nonumber\\
	& l_{1,\alpha}: \mathcal{Y}\to \mathbb{R},  
	\quad
	l_{1,\alpha}(y) := (1-\alpha)\jjntqT \rho^{-2}_0y_{obs} \, y\, dxdt, 
	\nonumber\\
	& l_{2,\alpha}: \Lambda\to \mathbb{R},  
	\quad 
	l_{2,\alpha}(\lambda) := -\alpha \jjntqT  y_{obs}\, L^\star (\rho^{-1}\lambda)\,dxdt.  \nonumber
\end{align}

	From the symmetry of $a_{r,\alpha}$ and $c_{\alpha}$, we easily check that this formulation corresponds 
	to the saddle point problem\,: 
\[
\left\{
	\begin{aligned}
	& \sup_{\lambda\in \Lambda} \inf_{y\in \mathcal{Y}} \mathcal{L}_{r,\alpha}(y,\lambda), \\
	& \mathcal{L}_{r,\alpha}(y,\lambda):=\mathcal{L}_r(y,\lambda) - \frac{\alpha}{2} \biggl\Vert \rho_0 \biggl(L^{\star}(\rho^{-1}\lambda)
	    +\rho_0^{-2}(y-y_{obs})1_{\om}\biggr)\biggr\Vert_{L^2(Q_T)}^2. 
	\end{aligned}
\right.
\]
\begin{proposition}\label{prop:mfalpha}
	Let $\rho_0\in \mathcal{R}$ and $\rho\in \mathcal{R}\cap L^{\infty}(Q_T)$. Then, for any $\alpha\in (0,1)$ and $r>0$, the formulation \eqref{eq:mfalpha} is well-posed. Moreover, the unique pair 
	$(y,\lambda)$ in $\mathcal{Y}\times \Lambda$ satisfies 
\begin{equation}\label{estimate_solalpha}
	\theta_1 \Vert y\Vert_\mathcal{Y}^2 + \theta_2 \Vert \lambda\Vert_{\Lambda}^2 
	\leq  \biggl(\frac{(1-\alpha)^2}{\theta_1} + \frac{\alpha^2}{\theta_2}\biggr)\Vert \rho_0^{-1}y_{obs}\Vert^2_{L^2(q_T)}
\end{equation}
	with 
\begin{equation}
	\theta_1:=\min\biggl(1-\alpha, \frac{r}{\eta}\biggr) , \quad \theta_2:= \frac{\alpha \rho_{\star}}{\rho_{\star} + C_{\Omega, T}}   \nonumber
\end{equation}
where $C_{\Omega,T}$ is the continuity constant so that $\Vert \rho^{-1}\lambda\Vert_{L^2(Q_T)}\leq C_{\Omega,T}\Vert L^{\star}(\rho^{-1}\lambda)\Vert_{L^2(Q_T)}$ for any $\lambda\in \Lambda$.
\end{proposition}
\textsc{Proof-} We easily get the continuity of the bilinear forms $a_{r,\alpha}$, $b_{\alpha}$ and $c_{\alpha}$: 
 $$
\begin{aligned}
& \vert a_{r,\alpha}(y,\overline{y})\vert 	\leq \max(1-\alpha,r\eta^{-1}) \Vert y\Vert_\mathcal{Y} \Vert\overline{y}\Vert_\mathcal{Y}, \quad \forall y,\overline{y}\in \mathcal{Y}, \\
& \vert b_{\alpha}(y,\lambda)\vert  \leq \max(\alpha,\eta^{-1/2}\|\rho\|_{L^{\infty}(Q_T)}) \Vert y\Vert_\mathcal{Y} \Vert \lambda\Vert_{\Lambda}, \quad \forall y\in \mathcal{Y},\,~\forall \lambda\in\Lambda  ,\\
& \vert c_{\alpha}(\lambda,\overline{\lambda})\vert \leq \alpha \Vert \lambda\Vert_{\Lambda} \Vert \overline{\lambda}\Vert_{\Lambda}, \quad \forall \lambda, \overline{\lambda}\in \Lambda 
\end{aligned}
$$
	and the continuity of the linear form $l_{1,\alpha}$ and $l_{2,\alpha}$~: 
$$
\Vert l_{1,\alpha}\Vert_{\mathcal{Y}^{\prime}}\le (1-\alpha)\Vert\rho_0^{-1} y_{obs}\Vert_{L^2(q_T)}\quad\hbox{and}
\quad\Vert l_{2,\alpha}\Vert_{\Lambda^{\prime}} \le \alpha\Vert\rho_0^{-1} y_{obs}\Vert_{L^2(q_T)}.$$   
 
 	Moreover, since $\alpha\in (0,1)$, we also obtain the coercivity of $a_{r,\alpha}$ and of $c_{\alpha}$: precisely, we check 
	that $a_{r,\alpha}(y,y) \geq \theta_1 \Vert y\Vert^2_{\mathcal{Y}}$ for all $y\in \mathcal{Y}$ while, for any 
	$m\in (0,1)$, by writing 
\[
	\begin{aligned}
		c_{\alpha}(\lambda,\lambda) = &~ \alpha \Vert \rho_0L^\star(\rho^{-1}\lambda)\Vert^2_{L^2(Q_T)} 
		 = \alpha m \|\rho_0L^\star(\rho^{-1}\lambda)\|_{L^2(Q_T)}^2 
		 + \alpha(1-m) \|\rho_0L^\star(\rho^{-1}\lambda)\|_{L^2(Q_T)}^2 \\
		 \noalign{\smallskip}\dis
		\ge &~ \alpha m \|\rho_0L^\star(\rho^{-1}\lambda)\|_{L^2(Q_T)}^2 
		+ \frac{\alpha(1 - m)\rho_{\star}}{C_{\Omega, T}} \|\rho^{-1}\lambda\|_{L^2(Q_T)}^2 
		\ge \alpha\min\left(m, \ {(1-m)\rho_{\star}\over C_{\Omega, T}}\right)\|\lambda\|_{\Lambda}^2,
	\end{aligned}
\]
	we get $c_{\alpha} (\lambda, \lambda) \geq \theta_2 \|\lambda\|_{\Lambda}^2$ for all $\lambda\in \Lambda$ with $m = \rho_{\star}(\rho_{\star} + C_{\Omega, T})^{-1}$.

The result \cite[Prop 4.3.1]{brezzi_new}  implies the well-posedness of the mixed formulation \eqref{eq:mfalpha} and the estimate (\ref{estimate_solalpha}). \Fin

The $\alpha$-term in $\mathcal{L}_{r,\alpha}$ is a stabilization term\,: it ensures a coercivity property of $\mathcal{L}_{r,\alpha}$ with respect to the variable $\lambda$ and automatically the well-posedness, assuming here $r>0$. In particular, there is no need to prove any inf-sup property for the application $b_{\alpha}$.

\begin{proposition}
	The solutions of \eqref{eq:mf} and \eqref{eq:mfalpha} coincide. 
\end{proposition}
\textsc{Proof-} From the optimality system (\ref{eq:system_lambda}),  the multiplier $\lambda$ such that $(y,\lambda)$ is the solution of (\ref{eq:mf}) belongs to the more regular space $\Lambda$.  Therefore,  this pair $(y,\lambda)$ is also a solution of  the mixed formulation \eqref{eq:mfalpha}. 
	The result then follows from the uniqueness of the two formulations. \Fin

\subsection{Dual formulation of the extremal problem (\ref{eq:mf})} \label{sec1_dual}

As discussed at length in \cite{NC-AM-mixedwave}, we may also associate to the extremal problem \eqref{P} a equivalent problem involving only the variable $\lambda$. Again, this is particularly interesting at the numerical level. This requires a strictly positive augmentation parameter $r$.

For any $r>0$, let us define the linear operator $\mathcal{T}_r$ from $L^2(Q_T)$ into $L^2(Q_T)$ by 
\[
\mathcal{T}_r\lambda:=  \rho^{-1}Ly, \quad \forall \lambda\in L^2(Q_T)
\]
where $y \in \mathcal{Y}$ is the unique solution to
\begin{equation}\label{eq:imageA}
a_r(y, \overline y) = b(\overline y, \lambda), \quad \forall \overline y \in \mathcal{Y}.   
\end{equation}
The assumption $r>0$ is necessary here in order to guarantee the well-posedness of (\ref{eq:imageA}). Precisely, for any $r>0$, the form $a_r$ defines a norm equivalent to the norm on $\mathcal{Y}$.

The following important lemma holds\,:
\begin{lemma}\label{lemmaA}
For any $r>0$, the operator $\mathcal{T}_r$ is a strongly elliptic, symmetric isomorphism from $L^2(Q_T)$ into $L^2(Q_T)$.
\end{lemma}
\textsc{Proof-} From the definition of $a_r$, we easily get that $\Vert \mathcal{T}_r\lambda\Vert_{L^2(Q_T)}\leq r^{-1} \Vert \lambda\Vert_{L^2(Q_T)}$ and the continuity of $\mathcal{T}_r$. Next, consider any $\lambda^{\prime}\in L^2(Q_T)$ and denote by $y^{\prime}$ the corresponding unique solution of (\ref{eq:imageA}) so that $\mathcal{T}_r\lambda^{\prime}:=\rho^{-1}Ly^{\prime}$. Relation (\ref{eq:imageA}) with $\overline{y}=y^{\prime}$ then implies that 
\begin{equation}
\jjntQT (\mathcal{T}_r\lambda^{\prime}) \lambda dx\,dt = a_r(y,y^{\prime})    \label{arAlambda}
\end{equation}
and therefore the symmetry and positivity of $\mathcal{T}_r$. The last relation with $\lambda^{\prime}=\lambda$ and the observability estimate (\ref{crucial_estimate_est}) imply that $\mathcal{T}_r$ is also positive definite.

Finally, let us check the strong ellipticity of $\mathcal{T}_r$, equivalently that the bilinear functional $(\lambda,\lambda^{\prime})\mapsto \jjntQT (\mathcal{T}_r\lambda)\lambda^{\prime} dx\,dt$ is $L^2(Q_T)$-elliptic. Thus, we want to show that 
\begin{equation}\label{ellipticity_A}
\jjntQT (\mathcal{T}_r\lambda)\lambda\,dx\,dt \geq C \Vert \lambda\Vert^2_{L^2(Q_T)}, \quad\forall \lambda\in L^2(Q_T)
\end{equation}
for some positive constant $C$. Suppose that (\ref{ellipticity_A}) does not hold; there exists then a sequence $\{\lambda_n\}_{n\geq 0}$ of $L^2(Q_T)$ such that 
\[
\Vert \lambda_n\Vert_{L^2(Q_T)}=1, \quad\forall n\geq 0, \qquad \lim_{n\to\infty} \jjntQT  (\mathcal{T}_r\lambda_n)\lambda_n \,dx\,dt=0.
\]
Let us denote by $y_n$ the solution of (\ref{eq:imageA}) corresponding to $\lambda_n$. From (\ref{arAlambda}), we then obtain that 
\begin{equation}
\lim_{n\to\infty} \left(r \Vert \rho^{-1}Ly_n\Vert^2_{L^2(Q_T)}+ \Vert \rho_0^{-1}y_n\Vert^2_{L^2(q_T)}\right)=0 .\label{limit}
\end{equation}
From (\ref{eq:imageA}) with $y=y_n$ and $\lambda=\lambda_n$, we have
\begin{equation} \label{phinlambdan}
\jjntQT   (r\rho^{-1}  L y_n - \lambda_n) \rho^{-1}L\overline{y} \,dx\,dt + \jjntqT \rho_0^{-2}y_n \overline{y} \,dx\,dt =0,
\end{equation}
for every $\overline{y} \in \mathcal{Y}$.
We define the sequence $\{\overline{y}_n\}_{n\geq 0}$ as follows : 
\begin{equation}\nonumber
\left\{
\begin{aligned}
& \rho^{-1}L\overline{y}_n = r\, \rho^{-1} Ly_n-\lambda_n, & &\textrm{in}\quad Q_T, \\
& \overline{y}_n=0, & &\textrm{in}\quad \Sigma_T,\\
& \overline{y}_n(\cdot,0)=0, & &\textrm{in}\quad \Omega,
\end{aligned}
\right.
\end{equation}
so that, for all $n$, $\overline{y}_n$ is the solution of the heat equation with zero initial data and source term $r \rho^{-1}Ly_n -\lambda_n$ in $L^2(Q_T)$. Energy estimates imply that 
$$\Vert \rho_0^{-1}\overline{y}_n\Vert_{L^2(q_T)} \leq C_{\Omega,T} \rho_{\star}^{-1}\Vert \rho\Vert_{L^{\infty}(Q_T)}  \Vert r \rho^{-1}Ly_n-\lambda_n\Vert_{L^2(Q_T)}$$ 
and that $\overline{y}_n\in \mathcal{Y}$. Then, using (\ref{phinlambdan}) with $\overline{y}=\overline{y}_n$ we get 
\[
\Vert r\rho^{-1}Ly_n-\lambda_n\Vert_{L^2(Q_T)} \leq C_{\Omega,T}\rho_{\star}^{-1}\Vert \rho\Vert_{L^{\infty}(Q_T)} \Vert \rho_0^{-1}y_n\Vert_{L^2(q_T)}.
\]
Then, from (\ref{limit}), we conclude that $\lim_{n\to +\infty} \Vert \lambda_n\Vert_{L^2(Q_T)}=0$ leading to a contradiction and to the strong ellipticity of the operator $\mathcal{T}_r$. \Fin

The introduction of the operator $\mathcal{T}_r$ is motivated by the following proposition\,: 

\begin{proposition}\label{prop_equiv_dual}
For any $r>0$, let $y_0\in \mathcal{Y}$ be the unique solution of 
\[
a_r(y_0,\overline{y})= l(\overline{y}), \quad \forall \overline{y}\in \mathcal{Y}
\]
and let $J_r^{\star\star}:L^2(Q_T)\to L^2(Q_T)$ be the functional defined by 
\[
J_r^{\star\star}(\lambda) := \frac{1}{2} \jjntQT (\mathcal{T}_r \lambda)\lambda \,dx\,dt - b(y_0, \lambda).
\]
The following equality holds : 
\begin{equation}\nonumber
\sup_{\lambda\in L^2(Q_T)}\inf_{y\in \mathcal{Y}} \mathcal{L}_r(y,\lambda) = - \inf_{\lambda\in L^2(Q_T)} J_r^{\star\star}(\lambda)\quad + \mathcal{L}_r(y_0,0).
\end{equation}
\end{proposition}
The proof is standard and we refer for instance to \cite{NC-AM-mixedwave} in a similar context.
This proposition reduces the search of $y$, solution of problem \eqref{P}, to the minimization of $J_r^{\star\star}$ with respect to $\lambda$. This extremal problem is well-posed in view of Lemma \ref{lemmaA} and the ellipticity of the operator $\mathcal{T}_r$. 

%\textsc{Proof-} For any $\lambda\in L^2(Q_T)$, let us denote by $\ph_{\lambda}\in \Phi$ the minimizer of $\ph\to \mathcal{L}_r(\ph,\lambda)$; $\ph_{\lambda}$ satisfies the equation
%\[
%a_r(\ph_{\lambda},\overline{\ph}) + b(\overline{\ph},\lambda) = l(\overline{\ph}), \quad \forall \overline{\ph}\in \Phi
%\]
%and can be decomposed as follows : $\ph_{\lambda} = \psi_{\lambda} + \ph_0$ where $\psi_{\lambda}\in \Phi$ solves 
%\[
%a_r(\psi_{\lambda},\overline{\ph}) + b(\overline{\ph},\lambda) = 0, \quad \forall \overline{\ph}\in \Phi.
%\]
%We then have 
%\begin{equation}
%\nonumber
%\begin{aligned}
%\inf_{\ph\in \Phi} \mathcal{L}_r(\ph,\lambda)  & = \mathcal{L}_r(\ph_\lambda,\lambda) = \mathcal{L}_r(\psi_\lambda+ \ph_0,\lambda)  \\
%&= \frac{1}{2}a_r(\psi_\lambda+ \ph_0,\psi_\lambda+ \ph_0) + b(\psi_\lambda+ \ph_0,\lambda) - l(\psi_\lambda+ \ph_0) \\
%& := X_1 +X_2 + X_3
%\end{aligned}
%\end{equation}
%with
%%
%\begin{equation}
%\nonumber
%\left\{
%\begin{aligned}
%& X_1 = \frac{1}{2}a_r(\psi_\lambda,\psi_\lambda) + b(\psi_\lambda,\lambda) + b(\ph_0,\lambda)  \\
%& X_2 = a_r(\psi_\lambda,\ph_0) - l(\psi_\lambda), \quad X_3=  \frac{1}{2}a_r(\ph_0,\ph_0) - l(\ph_0).
%\end{aligned}
%\right.
%\end{equation}
%From the definition of $\ph_0$, $X_2=0$ while $X_3= \mathcal{L}_r(\ph_0,0)$. Eventually, from the definition of $\psi_\lambda$, 
%\[
%X_1 = -\frac{1}{2}a_r(\psi_{\lambda},\psi_{\lambda}) + b(\ph_0,\lambda) = -\frac{1}{2}\jjntQT (\mathcal{P}_r \lambda) \,\lambda\, dx \,dt + b(\ph_0,\lambda)
%\]
%and the result follows. \Fin

\begin{remark}
Assuming in addition that the domain $\Omega$ is of class $C^2$, the results of this section apply if the distributed observation on $q_T$ is replaced by a Neumann boundary observation on the open subset $\gamma$ of 
$\partial\Omega$ (i.e.  assuming $y_{obs}:=\frac{\partial y}{\partial \nu}\in L^2(\gamma_T)$ is known on $\gamma_T:=\gamma\times(0,T)$). This is due to the following Carleman inequality, proved in \cite{FursikovImanuvilov} : there exists a positive constant $C=C(\omega,T,\Vert c\Vert_{C^1(\overline{\Omega})},\Vert d\Vert_{L^{\infty}(Q_T)})$ such that 
\begin{equation}\label{boundary_carleman}
	\Vert\tilde\rho_{c,0}^{-1}y\Vert^2_{L^2(Q_T)}+\Vert \tilde\rho_{c,1}^{-1}\nabla y\Vert^2_{L^2(Q_T)}
	 \leq C \Vert y\Vert_{\tilde{\mathcal{Y}}_0}^2,  
\end{equation}
	for any $y\in   \widetilde{\mathcal{Y}}_0:=\left\{ y \in C^2(\overline Q_T): y = 0 \textrm{ on }  \Sigma_T\setminus\overline{\gamma}_T\right\}$, where
\[
	(y, \ybar)_{\widetilde{\mathcal{Y}}_0} = \int\!\!\!\!\int_{\gamma_T} \tilde\rho^{-2}_{c,1} \frac{\partial y}{\partial \nu} \,\frac{\partial \ybar}{\partial \nu}\, d\Gamma\,dt + \jjntQT \tilde\rho^{-2}_cL\, y L\, \ybar\, dx\,dt
\]
	and $\Vert y\Vert_{\tilde{\mathcal{Y}}_0}^2=(y, y)_{\widetilde{\mathcal{Y}}_0}$. Here, 
	$\tilde\rho_c$, $\tilde\rho_{c,0}$ and $\tilde\rho_{c,1}$ are appropriate weight functions similar
	to $\rho_c$, $\rho_{c,0}$ and $\rho_{c,1}$ respectively (see \eqref{weights}).
 
%\Vert y(\cdot,0),y_t(\cdot,0)\Vert^2_{H_0^1(\Omega)\times L^2(\Omega)} \leq C_{obs} \biggl( \biggl\Vert \frac{\partial y}{\partial %\nu}\biggr\Vert^2_{L^2(\Sigma_T)}+ \Vert Ly\Vert^2_{L^2(Q_T)}  \biggr), \quad \forall y\in \hat Z
%\end{equation}
%\forall (y,\pvec) \in   \mathcal{B}_{\rho_{p,0},\rho_{p,1},\rho_p}
	Actually, it suffices to re-define the forms $a$ and $l$ in \eqref{eq:mf} by 
	$$
	\tilde a(y,\overline{y}):=\int\!\!\!\!\int_{\gamma_T} \tilde\rho^{-2}_{c,1} \frac{\partial y}{\partial \nu} \,\frac{\partial \ybar}{\partial \nu}\, d\Gamma\,dt
	\quad\hbox{and}\quad 
	\tilde l(y):=\int\!\!\!\!\int_{\gamma_T} \tilde\rho^{-2}_{c,1} \frac{\partial y}{\partial \nu} \, y_{obs}\, d\Gamma\,dt \quad \forall
	y,\ybar\in \widetilde{\mathcal{Y}},$$
	where $\widetilde{\mathcal{Y}}$ is the completion of  $\widetilde{\mathcal{Y}}_0$ with respect to 
	the scalar product $(\cdot,\cdot)_{\tilde{\mathcal{Y}}_0}$.
\end{remark}

\begin{remark}\label{rk_control_analogy}
We also emphasize that the mixed formulation (\ref{eq:mf}) has a structure very close to the one we get when we address - using the same approach - the null controllability of (\ref{eq:heat}): more precisely, the control of minimal $L^2(q_T)$-norm which drives to rest the initial data $y_0\in L^2(\Omega)$ is given by $v=\rho_0^{-2}\ph\, 1_{q_T}$ where $(\ph,\lambda)\in \Phi\times L^2(Q_T)$ solves the mixed formulation
\begin{equation} \nonumber %\label{eq:mfi}
\left\{
\begin{array}{rcll}
\noalign{\smallskip} a(\ph, \overline{\ph}) + b(\overline{\ph}, \lambda) & = & l(\overline{\ph}), & \qquad \forall \overline{\ph} \in \Phi,\\
\noalign{\smallskip} b(\ph, \overline{\lambda}) & = & 0, & \qquad \forall \overline{\lambda} \in L^2(Q_T),
\end{array}
\right.
\end{equation}
where
\begin{align}
& a : \Phi \times \Phi \to \mathbb{R},  \quad a(\ph, \overline{\ph}) := \jjntqT \rho_0^{-2}\ph(x, t) \overline{\ph}(x, t) \, dx \, dt, \nonumber\\
& b: \Phi \times L^2(Q_T)  \to \mathbb{R},  \quad b(\ph, \lambda) := \jjntQT \rho^{-1}L^{\star}\ph\,\lambda\,dx dt, \nonumber\\
& l: \Phi \to \mathbb{R},  \quad l(\ph) := -(\ph(\cdot, 0),y_0)_{L^2(\Omega)}. \nonumber
\end{align}
Here, the weights $\rho$ and $\rho_0$ are taken in a space of functions that may blow-up at time $t=T$ and
$\Phi$ is a complete space associated to these weights. For more details, see \cite{munch_desouza}.
%with here $\overline{\rho}(x,t)=\rho(x,T-t)$ and $\overline{\rho}_0(x,t)=\rho_0(x,T-t)$. We refer to \cite{munch_desouza}.
\end{remark}

\begin{remark}
Reversing the order of priority between the constraints  $y-y_{obs}=0$ in $L^2(q_T)$ and $\rho^{-1}(Ly-f)=0$ in $L^2(Q_T)$, a possibility could be to minimize the functional 
$y\to \Vert \rho^{-1}(Ly-f)\Vert_{L^2(Q_T)}$ over $y\in\mathcal{Y}$ subject to the constraint $\rho_0^{-1}(y-y_{obs})=0$ in $L^2(q_T)$ via the introduction of a Lagrange multiplier in $L^2(q_T)$. 
The fact that the following inf-sup property : there exists $\delta>0$ such that 
$$
\inf_{\lambda\in L^2(q_T)} \sup_{y\in\mathcal{Y}} \frac{\jjntqT \rho^{-1}_0 y \lambda \,dxdt }{\Vert \lambda\Vert_{L^2(q_T)} \Vert y\Vert_{\mathcal{Y}_0}}   \geq \delta
$$
associated to the corresponding mixed-formulation holds true is however an open issue. 
On the other hand, if a $\eps$-term is added as in (\ref{FV}), this property is satisfied (we refer again to the book \cite{Klibanov-book}).
\end{remark}

%\begin{remark}
%{\color{blue} find in the literature something in order to address the reconstruction of a source term as well } 
%\end{remark}

%%%% NEW SECTION

%%%%%%%%%%%%%%%%%
%%%%%%%%%%%%%%%%%
\section[First order mixed formulation]{Recovering the solution from a partial observation: a first order mixed re-formulation}\label{recovering_y_fo}
%%%%%%%%%%%%%%%%%
%%%%%%%%%%%%%%%%%

In this section, we consider a first order mixed formulation of the parabolic equation (\ref{eq:heat}) introducing the flux variable $\pvec:=c(x)\nabla y$. We then apply to this first order system the methods developed in the previous section and address the reconstruction of $y$ and $\pvec$ from the distributed observation $y_{obs}$. The introduction of this equivalent first order system is advantageous at the numerical level as it allows to reduce the regularity order of the spaces. 

\subsection{Direct approach: Minimal local weighted $L^2$-norm}

 We rewrite the parabolic equation  (\ref{eq:heat}) as the following equivalent first order system : 
\begin{equation}
\label{eq:heat_mixed}
	\left\{
		\begin{array}{lll}
			y_t- \nabla\cdot\pvec+d\,y  =f, \quad c(x)\nabla y-\pvec=0 	& \textrm{in}& Q_T,    		\\
   			y  = 0										  		&\textrm{on}& \Sigma_T, 	\\
   			y(x, 0) = y_0(x) 								  		& \textrm{in}& \Om.
   		\end{array} 
 	\right.
\end{equation}
	
The reformulation of the parabolic equation (\ref{eq:heat}) into a first order system is standard and has been analyzed for instance in \cite{jerome_space_time,li_todd_huang}: there, the existence and uniqueness of solution for a associated $L^2-H(div)$ weak formulation is proved, together with a priori estimates assuming notably that $y_0\in H_0^1(\Omega)$. In the Appendix, assuming only $y_0\in L^2(\Omega)$, we study the well-posedness of a $H_0^1-L^2$ weak formulation associated to the problem (\ref{eq:heat_mixed}). We refer to Proposition \ref{ex_un_mixed_sol1}.  	

%In \cite{jerome_space_time}, it is for instance prove that if $f\in L^2(Q_T)$ and $y_0\in H_0^1(\Omega)$, then problem (\ref{eq:heat_mixed}) has a unique solution such that
%%
%\begin{equation}
%(y,\pvec)\in H^1(0,T; L^2(\Omega))\times (L^2(0,T; H(div,\Omega))\cap L^{\infty}(0,T; \mathbf{L^2(\Omega))}).
%\end{equation}
%%	
In the sequel, we use the following notations : 
\begin{equation}
\mathcal{I}(y,\pvec):=y_t- \nabla\cdot\pvec+d\,y, \qquad \mathcal{J}(y,\pvec):=c(x)\nabla y-\pvec.
\end{equation}

Assuming again for simplicity that $f=0$ and proceeding as in the previous, we address the reconstruction of $y$ and now $\pvec$ from the observation $y_{obs}$ by introducing a least-squares type problem. 

Precisely, we first define the space 
\[
	 \mathcal{U}_0 = \left\{ (y,\pvec)\in C^1(\overline Q_T)\times \Cvec^1(\overline Q_T): y = 0 \textrm{ on }  \Sigma_T\right\}
\]
and for any $\eta_1,\eta_2>0$ and any $\rho,\,\rho_0,\,\rho_1\in \mathcal{R}$, we define the bilinear form 
\[
\begin{alignedat}{2}
	((y,\pvec), (\ybar,\overline{\pvec}))_{\mathcal{U}_0} =&~ \jjntqT\rho_0^{-2} y \,\ybar\, dx\,dt 
	+ \eta_1\jjntQT\rho^{-2}_1 \mathcal{J}(y,\pvec)\cdot \mathcal{J}(\ybar,\overline{\pvec})\, dx\,dt\\
	&~+ \eta_2\jjntQT\rho^{-2}\mathcal{I}(y,\pvec)\mathcal{I}(\ybar,\overline\pvec)\, dx\,dt\quad\forall (y,\pvec),(\ybar,\overline{\pvec}) \in \mathcal{U}_0.
\end{alignedat}
\]

	From unique continuation properties for parabolic equations, this bilinear form defines a scalar product (we also refer to Proposition \ref{carleman_puel} which quantifies the unique continuation). 
	Let then $ \mathcal{U}$ be the completion of $\mathcal{U}_0$ for this scalar 
	product and denote the norm over $ \mathcal{U}$ by  $\Vert\cdot\Vert_\mathcal{U}$ 
	such that
\begin{equation}
	\Vert (y,\pvec)\Vert^2_\mathcal{U}:=\Vert \rho_0^{-1}y\Vert^2_{L^2(q_T)}
	+ \eta_1 \Vert \rho_1^{-1} \mathcal{J}(y,\pvec) \Vert^2_{\Lvec^2(Q_T)} 
	+ \eta_2 \Vert \rho^{-1}\mathcal{I}(y,\pvec) \Vert^2_{L^2(Q_T)}.
\end{equation}

	Finally, we define the closed subset $\mathcal{V}$ of $\mathcal{U}$ by 
\begin{equation}
	\mathcal{V}:=\biggl\{ (y,\pvec)\in \mathcal{U}\,:\,
	\rho^{-1}_1\mathcal{J}(y,\pvec)=0\; \textrm{in}\; \Lvec^2(Q_T)\;\hbox{ and }\;\rho^{-1}\mathcal{I}(y,\pvec)=0\; \textrm{in}\; L^2(Q_T)\biggr\}
\end{equation}
	and we endow $\mathcal{V}$ with the same norm than $\mathcal{U}$.

The following problem into consideration is then as follows\,: 
\begin{equation}
\label{inverse_pb_2}
	\left\{
		\begin{array}{l}
			\dis \hbox{Minimize }\ J(y,\pvec) := {1 \over 2} \jjntqT \rho^{-2}_0 |y(x,t)-y_{obs}(x,t)|^2 \,dx\,dt \\
			\noalign{\smallskip}
			\hbox{Subject to }\ (y,\pvec) \in \mathcal{V}.
		\end{array}
	\right.
\end{equation}
As in Section \ref{sec1_direct}, this extremal is well-posed in view of the definition of $\mathcal{V}$: there exists a unique pair $(y,\pvec)$, minimizer for $J$.

	The pair $(y,\pvec)$ of this extremal problem, equivalent to the previous one, are now submitted to the constraints 
	$\rho^{-1}_1\mathcal{J}(y,\pvec)=0$ (in $\Lvec^2(Q_T)$) and $ \rho^{-1}\mathcal{I}(y,\pvec)=0$ (in $L^2(Q_T)$). 
	As before, these constraints are addressed by introducing Lagrange multipliers. 
	
	Precisely, we set
	$\mathcal{X}:=L^2(Q_T)\times\Lvec^2(Q_T)$ and then we consider the following mixed
	formulation~: find $((y,\pvec),(\lambda,\boldsymbol{\mu}))\in\mathcal{U}\times\mathcal{X}$
	solution of 
\begin{equation} 
\label{eq:mf4}
	\left\{
		\begin{array}{rcll}
			\noalign{\smallskip} 
			a((y,\pvec), (\ybar,\overline{\pvec})) +  b((\ybar,\overline{\pvec}),(\lambda,\boldsymbol{\mu})) & = & 
			l(\ybar,\overline{\pvec}) & \quad \forall (\ybar,\overline{\pvec}) \in   \mathcal{U},\\
			\noalign{\smallskip}
			b( (y,\pvec), (\lambar,\mubar)) & = & 0 & \quad \forall (\lambar,\mubar) \in \mathcal{X},
		\end{array}
	\right.
\end{equation}
	where
\begin{align}
	\nonumber
	& a : \mathcal{U}\times \mathcal{U}\to \mathbb{R},  
	\quad  a((y,\pvec),(\ybar,\pvec))  :=   \jjntqT\rho_0^{-2} y \,\ybar\, dx\,dt,\\
	\nonumber
	& b: \mathcal{U} \times\mathcal{X}\to \mathbb{R},\quad b((y,\pvec), (\lambda,\boldsymbol{\mu})) := 
	\jjntQT\rho^{-1}_1 \mathcal{J}(y,\pvec)\cdot\boldsymbol{\mu}\,dx\,dt+\jjntQT\rho^{-1}\mathcal{I}(y,\pvec) \lambda\,dx\,dt\nonumber\\
	\nonumber
	& l: \mathcal{U}\to \mathbb{R}, 
	\quad  l(y,\pvec) :=\jjntqT \rho^{-2}_0y\, y_{obs}\, dx\,dt.
\end{align}

We have the following result\,: 
\begin{theorem}\label{th:mf4}
	Let $\rho_0\in \mathcal{R}$ and $\rho,\rho_1\in \mathcal{R}\cap L^{\infty}(Q_T)$. We have\,:
\begin{enumerate}
	\item The mixed formulation \eqref{eq:mf4} is well-posed.
	\item The unique solution $((y,\pvec),(\lambda,\boldsymbol{\mu}))\in\mathcal{U}\times\mathcal{X}$ is the unique saddle-point of 
	the Lagrangian ${\mathcal{L} :  \mathcal{U}\times \mathcal{X}\to\mathbb{R}}$ defined by
\begin{equation}
\label{eq:calLm4}
	\begin{alignedat}{2}
		\mathcal{L}((y,\pvec), (\lambda,\boldsymbol{\mu})) :=&~ {1\over 2} a((y,\pvec),(y,\pvec)) 
		+ b((y,\pvec),(\lambda,\boldsymbol{\mu})) -l(y,\pvec).
	\end{alignedat}
\end{equation}
	\item The unique solution $((y,\pvec),(\lambda,\boldsymbol{\mu}))$ satisfies the following estimate : 
\begin{equation}
\label{estimative_sol_weig_fo}
\begin{alignedat}{2}
	&\hspace{-0.8cm}\|(y,\pvec)\|_\mathcal{U}\leq\|\rho_0^{-1}y_{obs}\|_{L^2(q_T)},\\
	&\hspace{-0.8cm}\|(\lambda,\boldsymbol{\mu})\|_{\mathcal{X}}\leq2 \sqrt{\max\{C_{\Omega,T}\rho_{\star}^{-2}\Vert \rho_1\Vert_{L^{\infty}(Q_T)}^2+\eta_1,C_{\Omega,T}\rho_{\star}^{-2}\Vert \rho\Vert_{L^{\infty}(Q_T)}^2+\eta_2\}}\,\|\rho_0^{-1}y_{obs}\|_{L^2(q_T)}.
	\end{alignedat}
\end{equation}
\end{enumerate}
\end{theorem}
\textsc{Proof-} 
	The proof is similar to the proof of Theorem \ref{th:mf}. From the definition, the bilinear form $ a$ is continuous over 
	$ \mathcal{U}\times \mathcal{U}$, symmetric and positive and the bilinear form $ b$ is continuous 
	over $\mathcal{U}\times\mathcal{X}$. Furthermore, the linear form $ l$ is continuous over 
	$ \mathcal{X}$. In particular, we get
\begin{equation}\label{norm_linear_form_weig_yp}
	\|l\|_{ \mathcal{X}'}\leq\|\rho^{-1}_0y_{obs}\|_{L^2(q_T)},\quad
	 \|a\|_{\mathscr{L}^2(\mathcal{U})}\leq1,\quad
	 \|b\|_{\mathscr{L}^2(\mathcal{U},\mathcal{X})}\leq\max\{\eta_1^{-{1/2}},\eta_2^{-{1/2}}\}.
\end{equation}

	Therefore, the well-posedness of the formulation \eqref{eq:mf4} is the consequence of two properties\,: 
	first, the coercivity of the form $ a$ on the kernel $\mathcal{N}( b):=\{(y,\pvec)\in \mathcal{U}\,: 
	b((y,\pvec),(\lambda,\boldsymbol{\mu}))=0\,~ \forall (\lambda,\boldsymbol{\mu})\in \mathcal{X}\}$. Again, this 
	holds true since the kernel of $b$ coincides with the space $\mathcal{V}$. 
	
	Second, the inf-sup property 
	which reads as\,: 
\begin{equation}
\label{eq:infsup4}
	\inf_{(\lambda,\boldsymbol{\mu}) \in \mathcal{X}} \sup_{(y,\pvec)\in \mathcal{U}} 
	\frac{ b((y,\pvec), (\lambda,\boldsymbol{\mu}))}{\|(y,\pvec)\|_\mathcal{U} \|(\lambda,\boldsymbol{\mu})\|_{\mathcal{X}}} \geq \delta
\end{equation}
for some $\delta >0$.

	Let us check this property. For any fixed $(\lambda^0,\boldsymbol{\mu}^0)\in \mathcal{X}$, we 
	define the (unique) element $(y^0,\pvec^0)$ such that 
\[
	\rho^{-1}\mathcal{I}(y^0,\pvec^0)=\lambda^0\;\textrm{ in } Q_T,
	\quad
	\rho_1^{-1}\mathcal{J}(y^0,\pvec^0)=\boldsymbol{\mu}^0\;\textrm{ in } Q_T, 
	\quad 
	y^0=0\textrm{ on } \Sigma_T, 
	\quad
	 y^0(\cdot,0)=0\;\textrm{ in } \Omega.
\]
	The pair $(y^0,\pvec^0)$ is therefore solution of a parabolic equation in the mixed form with source term 
	$(\rho\lambda^0,\rho_1\boldsymbol{\mu}^0)$ in $L^2(0,T;L^2(\Om))\times L^2(0,T;\Lvec^2(\Om))$,
	 null Dirichlet boundary condition and null initial state. From Proposition \ref{ex_un_mixed_sol1} applied with $f=\rho\lambda^0\in L^2(Q_T)$ and $\Fvec=\rho_1\boldsymbol{\mu}^0\in \Lvec^2(Q_T)$, the weak solution satisfies 
	 $(y^0,\pvec^0)\in (L^2(0,T;H^1_0(\Om))\cap C^0([0,T];L^2(\Om)))\times \Lvec^2(Q_T)$. Moreover, from (\ref{energy_est_mixed}), there exists 
	 a constant $C_{\Omega,T}$ such that the unique pair $(y^0,\pvec^0)$ satisfies the inequality
\begin{equation}
	\begin{alignedat}{2}\label{estimate_L2L2_4}
		\Vert \rho_0^{-1}y^0\Vert_{L^2(q_T)}^2\leq
		&~ C_{\Omega,T}\rho_{\star}^{-2}(\Vert \rho\Vert^2_{L^{\infty}(Q_T)}\|\lambda^0\|^2_{L^2(Q_T)}
		+\Vert \rho_1\Vert^2_{L^{\infty}(Q_T)} \|\boldsymbol{\mu}^0\|^2_{\Lvec^2(Q_T)})
	\end{alignedat}
\end{equation}	
	which proves that $(y^0,\pvec^0)\in \mathcal{U}$.

%	 and, for a.e. $t\in(0,T)$, it satisfies
%\[
%	\begin{array}{rcll}
%	\noalign{\smallskip} \dis
%	{d\over dt}(y(\cdot,t),w)- (\nabla\cdot\pvec(\cdot,t),w) + (d(\cdot,t)y(\cdot,t),w)
%	& = &(\rho(\cdot,t)\lambda^0(\cdot,t),w)  & \quad \forall w\in   L^2(\Omega),\\
%	\noalign{\smallskip} \dis
%	-(y(\cdot,t),\nabla\cdot\uvec)- (c^{-1}\pvec(\cdot,t),\uvec) & = &(\rho_1(\cdot,t)c^{-1}\mu^0(\cdot,t),\uvec)
%	& \quad \forall \uvec\in H(div,\Omega).
%		\end{array}
%\]	 
	
Consequently, we may write that  
\[
	\begin{alignedat}{2}
	\sup_{(y,\pvec)\in \mathcal{U}} \frac{b((y,\pvec), (\lambda^0,\boldsymbol{\mu}^0))}
	{\|(y,\pvec)\|_\mathcal{U} \|(\lambda^0,\boldsymbol{\mu}^0)\|_{\mathcal{X}}} 
	\geq\frac{b((y^0,\pvec^0), (\lambda^0,\boldsymbol{\mu}^0))}
	{\|(y^0,\pvec^0)\|_\mathcal{U} \|(\lambda^0,\boldsymbol{\mu}^0)\|_{\mathcal{X}}}\\
	= \frac{\|(\lambda^0,\boldsymbol{\mu}^0)\|_{\mathcal{X}}}{\left( \|\rho^{-1}_0y^0\|^2_{L^2(q_T)} +
	\eta_1 \|\boldsymbol{\mu}^0\|^2_{\Lvec^2(Q_T)}+ \eta_2 \|\lambda^0\|^2_{L^2(Q_T)} \right)^{1/2}}.
	\end{alignedat}
\]
leading together with (\ref{estimate_L2L2_4}) to 
\[
	\sup_{(y,\pvec)\in \mathcal{U}} \frac{b((y,\pvec), (\lambda^0,\boldsymbol{\mu}^0))}
	{\|(y,\pvec)\|_\mathcal{U} \|(\lambda^0,\boldsymbol{\mu}^0)\|_{\mathcal{X}}} 	\geq \delta,
\]
	with $\delta: = \left(\max\{C_{\Omega,T}\rho_{\star}^{-2}\Vert \rho_1\Vert_{L^{\infty}(Q_T)}^2+\eta_1,C_{\Omega,T}\rho_{\star}^{-2}\Vert \rho\Vert_{L^{\infty}(Q_T)}^2+\eta_2\}\right)^{-1/2}$. Hence, \eqref{eq:infsup4} holds.

	The point $(ii)$ is again due to the positivity and symmetry of the form $ a$.
	
	The point $(iii)$ is a consequence of classical estimates (see \cite{brezzi_new}, Theorem 4.2.3)\,:
\[
	\|(y,\pvec)\|_\mathcal{U}\leq{1\over\alpha_0}\|l\|_{\mathcal{U}'},
	\quad\|(\lambda,\boldsymbol{\mu})\|_{\mathcal{X}}\leq{1\over\delta}\left(1+{\|a\|_{\mathscr{L}^2(\mathcal{U})}\over\alpha_0}\right)
	\|l\|_{\mathcal{U}'},
\]
	where
\[
	\alpha_0:=\inf\limits_{y\in \mathcal{N}(b)} {a((y,\pvec),(y,\pvec))\over\|(y,\pvec\|_\mathcal{U}^2}.
\]	
	Estimates \eqref{norm_linear_form_weig_yp} and the equality $\alpha_0=1$ lead to the results.\Fin
	
	Again, as in Section \ref{sec1_direct}, the solution of (\ref{eq:mf4}) does not depend on the parameters $\eta_1,\eta_2$, only introduced in order to construct a scalar product in $\mathcal{U}_0$. In particular, $\eta_1$ and $\eta_2$ can be arbitrarily small. 

	Now, let us recall the following important result, analogue of Proposition \ref{carleman_estimate}, which provides a global estimate of $y$, solution of a parabolic equation with $L^2(H^{-1})$ right hand side, from a local  (in $q_T$) observation. 	% proved in Theorem $2.2$ in \cite{imanuv_puel_yama_2} 
	%and in Theorem $2.1$ in \cite{imanuv_puel_yama}.  
\begin{proposition}[Theorem $2.2$ in \cite{imanuv_puel_yama_2}]\label{carleman_puel}
	Let the weights $\rho_p,\rho_{p,0},\rho_{p,1}  \in \mathcal{R}$ $($see \eqref{def_spaceR}$)$ be defined  as follows\,:
 \begin{equation}
\label{weights_puel}
\begin{aligned}
 &\rho_p(x,t):=  \exp\left({\beta(x) \over t^2}\right), \quad \beta(x):=  K_{1}\left(e^{K_{2}} - e^{\beta_{0}(x)}\right), \\
 &\rho_{p,0}(x,t):=  t\rho_p(x,t), \quad \rho_{p,1}(x,t):=  t^{-1}\rho_p(x,t) ,\quad\rho_{p,2}(x,t):=  t^{-2}\rho_p(x,t) 
\end{aligned}
\end{equation}
 with $\beta_{0} \in C^\infty(\overline{\Omega})$ and where the positive constants $K_{i}$ are sufficiently large (depending on $T$, $c_0$, $\|c\|_{C^1(\overline{\Omega})}$ and $\|d\|_{L^\infty(Q_T)}$$)$ such that 
\[
\beta > 0\quad\textrm{in}\quad \Omega,\, \beta=0\quad\textrm{on}\quad \partial\Omega, \quad \nabla\beta(x)\neq0\quad\forall x\in\overline{\Omega}\setminus\omega.
\]
	Then, there exists a constant $C>0$, depending only on $\omega, T$, such that the following inequality holds\,: 
\begin{equation}%\label{crucial_estimate_h-1}
	\Vert  \rho_{p,0}^{-1}y\Vert^2_{L^2(Q_T)}+\Vert  \rho_{p,1}^{-1}\nabla y\Vert^2_{L^2(Q_T)}
		\leq C\left(\Vert\rho_p^{-1}\mathbf{G} \Vert_{\Lvec^2(Q_T)}^2+\Vert \rho_{p,2}^{-1}g \Vert_{L^2(Q_T)}^2
		+\Vert \rho_{p,0}^{-1}y\Vert^2_{L^2(q_T)}\right), \label{crucial_estimate_h}
\end{equation}
	where $y$ belongs to $\mathcal{K} := \left\{ y \in L^2(0,T;H_0^1 (\Om)): y_t \in L^2(0,T;H^{-1} (\Om))\right\}$ 
	and satisfies $Ly=g+\nabla\cdot\mathbf{G}$ in $Q_T$, with $g\in L^2(Q_T)$ and $\mathbf{G}\in \Lvec^2(Q_T)$.
\end{proposition}
This proposition allows to get the following second global estimate. 
	\begin{proposition}
	Let $\rho_p,\rho_{p,0},\rho_{p,1} \in \mathcal{R}$ the weights defined by (\ref{weights_puel}).
	There exists a constant $C>0$, depending only on $\omega,\,\Omega,\, T$,  such that 
\begin{equation} \label{crucial_estimate_puel}
	\Vert\rho_{p,0}^{-1}y\Vert^2_{L^2(Q_T)}+\Vert \rho_{p,1}^{-1}\nabla y\Vert^2_{L^2(Q_T)}
	+\Vert \rho_{p,1}^{-1}\pvec\Vert^2_{L^2(Q_T)} \leq C \Vert (y,\pvec) \Vert_{\mathcal{U}_p}^2
	 \quad\forall (y,\pvec) \in   \mathcal{U}_p,  
\end{equation}
	where $\mathcal{U}_p$ is the completion of $\mathcal{U}_0$  with respect to the scalar product
\[
	((y,\pvec), (\ybar,\overline{\pvec}))_{\mathcal{U}_{0,p}} \!= \!\!\!\jjntqT\!\!\!\rho^{-2}_{p,0} y \,\ybar\, dx\,dt
	 + \eta_1\!\jjntQT\!\!\!\!\rho^{-2}_{p,1} \mathcal{J}(y,\pvec)\cdot\mathcal{J}(\ybar,\overline{\pvec})\, dx\,dt
	+ \eta_2\!\jjntQT\!\!\!\!\rho^{-2}_p\mathcal{I}(y,\pvec)\mathcal{I}(\ybar,\overline\pvec)\, dx\,dt.
\]
\end{proposition}
\textsc{Proof-} First, let us to prove this inequality for $(y,\pvec) \in \mathcal{U}_0$. 
	So, let us introduce $\mathbf{G}:=\mathcal{J}(y,\pvec)$ and $g:=\mathcal{I}(y,\pvec)$. 
	Then, we also have that $y\in \mathcal{K}$ and \,$L\,y= g-\nabla\cdot\mathbf{G}$\, in \,$Q_T$.
	
	So, applying the Carleman inequality (\ref{crucial_estimate_h}), we obtain\,:
\[
	\begin{alignedat}{2}
		\Vert  \rho_{p,0}^{-1}y\Vert^2_{L^2(Q_T)}+\Vert  \rho_{p,1}^{-1}\nabla y\Vert^2_{L^2(Q_T)}
		\leq C\left(\Vert\rho_p^{-1}\mathbf{G} \Vert_{\Lvec^2(Q_T)}^2+\Vert \rho_{p,2}^{-1}g \Vert_{L^2(Q_T)}^2
		+\Vert \rho_{p,0}^{-1}y\Vert^2_{L^2(q_T)}\right).
	\end{alignedat}
\]
Moreover, writing that $\pvec=c(x)\nabla y-\Gvec$, we get 
\[
\begin{alignedat}{2}
	\Vert \rho_{p,1}^{-1}\pvec\Vert^2_{L^2(Q_T)} \leq& 2(\Vert \rho_{p,1}^{-1}c\nabla y\Vert^2_{L^2(Q_T)}
	+\Vert \rho_{p,1}^{-1}\mathbf{G}\Vert^2_{\Lvec^2(Q_T)}).
\label{crucial_estimate_h-1}
\end{alignedat}
\]	
Finally, since $\rho_{p,1}^{-1}\leq T\rho_p^{-1}$, we combine the last two inequalities to obtain
\[
	\Vert\rho_{p,0}^{-1}y\Vert^2_{L^2(Q_T)}+\Vert \rho_{p,1}^{-1}\nabla y\Vert^2_{L^2(Q_T)}
	+\Vert \rho_{p,1}^{-1}\pvec\Vert^2_{L^2(Q_T)} \leq C \Vert (y,\pvec) \Vert_{\mathcal{U}_0}^2
	 \quad\forall (y,\pvec) \in   \mathcal{U}_0.  
\]
	Then, by a density argument, we can deduce \eqref{crucial_estimate_puel}.\Fin
	
	\
	
Finally, assuming that the weights $\rho_0,\rho_1,\rho$, which appear in the mixed formulation (\ref{eq:mf4}), are related to the Carleman weights $\rho_{p,0},\rho_{p,1},\rho_p$ so that $\mathcal{U}\subset \mathcal{U}_p$,  we get the following stability result and a global estimate - analogue to (\ref{crucial_estimate_good}) - of any pair $(y,\pvec)\in \mathcal{U}$ in term of the norm $\Vert(y,\pvec)\Vert_{\mathcal{U}}$ (in particular for the solution of (\ref{eq:mf4})).
\begin{corollary}\label{stab_yp}
	Let $\rho_0\in \mathcal{R}$ and $\rho,\rho_1\in \mathcal{R}\cap L^{\infty}(Q_T)$ and
	assume that there exists a constant $K>0$ such that 
\begin{equation}
\rho_0\leq K \rho_{p,0}, \quad\rho_1\leq K \rho_{p,1}, \quad \rho\leq K \rho_{p,2} \quad \textrm{in}\quad Q_T.   \label{hypKyp}
\end{equation}
If $((y,\pvec),(\lambda,\mu))\in \mathcal{U}\times\mathcal{X}$ is the solution of the mixed formulation (\ref{eq:mf4}), then there exists $C>0$ such that 
\begin{equation}
\Vert \rho_{p,0}^{-1}y\Vert_{L^2(Q_T)}+\Vert \rho_{p,1}^{-1}\pvec\Vert_{L^2(Q_T)} \leq C \Vert( y,\pvec)\Vert_\mathcal{U}.  
 \label{crucial_estimate_good_fo}
 \end{equation}
\end{corollary}
\par\noindent
\textsc{Proof-} The hypothesis \eqref{hypKyp} implies that 
	$\mathcal{U}\subset \mathcal{U}_p$. 
	Therefore, estimate \eqref{estimative_sol_weig} and \eqref{crucial_estimate_puel} imply that 
\[
	\Vert \rho_{p,0}^{-1}y\Vert_{L^2(Q_T)}+\Vert \rho_{p,1}^{-1}\pvec\Vert_{L^2(Q_T)} 
	\leq C \Vert (y,\pvec)\Vert_{\mathcal{U}_p}
	\leq C \Vert( y,\pvec)\Vert_\mathcal{U}\leq C \Vert\rho_0^{-1} y_{obs}\Vert_{L^2(q_T)}.  
 \]
\Fin

Again, the functions $\rho_{p,0}^{-1}, \rho_{p,1}^{-1}$ vanish at time $t=0$, so that the variable $y$ and the flux $\pvec$ are reconstructed from the observation $y_{obs}$ everywhere in $Q_T$ except on the set $\Omega\times \{0\}$. Similarly, the weights $\rho_1$ and $\rho$ are introduced in the definition of $\mathcal{V}$ in order to reduce the effect of the "singularity" of the variable $y$ and $\pvec$ in the neighborhood of $\Omega\times \{t=0\}$. We refer to the discussion at the end of Section \ref{sec1_direct}. 

\begin{remark}\label{rk6}
As in Section \ref{recovering_y}, it is convenient to augment the lagrangien $\mathcal{L}$ defined in (\ref{eq:calLm4}) as follow: 
\begin{equation}
\nonumber
\left\{
\begin{aligned}
& \mathcal{L}_{\boldsymbol{r}}((y,\pvec),(\lambda,\boldsymbol{\mu})):=\frac{1}{2}a_r((y,\pvec),(y,\pvec))+b((y,\pvec),(\lambda,\boldsymbol{\mu}))-l(y,\pvec), \\
& a_{\boldsymbol{r}}((y,\pvec),(y,\pvec)):=a((y,\pvec),(y,\pvec))+ r_1 \Vert \rho_1^{-1}\mathcal{J}(y,\pvec)\Vert^2_{\Lvec^2(Q_T)}+ r_2 \Vert \rho^{-1}\mathcal{I}(y,\pvec)\Vert^2_{L^2(Q_T)}.
\end{aligned}
\right.
\end{equation}
for any $\boldsymbol{r}=(r_1,r_2)\in (\mathbb{R}^+)^2$.
\end{remark}

\begin{remark}\label{weak_lambdaMF}
Similarly to Remark \ref{weak_lambda},	the first equation of the mixed formulation \eqref{eq:mf4} reads as follows: 
\[
	\jjntqT\!\!\!\rho_0^{-2} y \,\ybar\, dx\,dt + \jjntQT\!\!\!\rho^{-1}_1 \mathcal{J}(\ybar,\overline{\pvec})\cdot\boldsymbol{\mu}\,dx\,dt
	+\jjntQT\!\!\!\rho^{-1}\mathcal{I}(\ybar,\overline{\pvec}) \lambda\,dx\,dt =\!\! \jjntqT\!\!\! \rho^{-2}_0\ybar\, y_{obs}\, dx\,dt
	\quad\forall (\ybar,\overline{\pvec}) \in \mathcal{U}.
\]

	But, according to Definition \ref{weak_heat_transposition}, this means that the pair $(\varphi,\boldsymbol{\sigma}):=(\rho^{-1}\lambda,c\rho_1^{-1}\boldsymbol{\mu})\in L^2(Q_T)\times \Lvec^2(Q_T)$ is solution 
	of the parabolic equation in the mixed form in the transposition sense, i.e.
	$(\varphi,\boldsymbol{\sigma})$ solves the problem\,: 
\begin{equation}
\label{eq:system_lambda_mu}
	\left\{
		\begin{array}{lcl}
   			\mathcal{I}^\star(\varphi,\boldsymbol{\sigma})=-\rho_0^{-2}\,(y-y_{obs})1_\om,
			\quad\mathcal{J}(\varphi,\boldsymbol{\sigma})=0  						&\textrm{in}&	Q_T,    	\\
			\varphi  = 0										&\textrm{on}&	\Sigma_T,	\\
			\varphi(\cdot,T)=0									&\textrm{in}&	\Om,
   		\end{array} 
 	\right.
\end{equation}
	where $\mathcal{I}^\star(\varphi,\boldsymbol{\sigma}):=-\varphi_t-\nabla\cdot\boldsymbol{\sigma}+d(x,t)\varphi$. In particular, this means that the multiplier pair $(\lambda,\boldsymbol{\mu})$, solution of a backward mixed system, vanishes if $y_{obs}$ is the restriction to $q_T$ of a solution of (\ref{eq:heat_mixed}). In this context, the rest of Remark \ref{weak_lambda}, in particular (\ref{QRbis}), can be adapted to (\ref{eq:mf4}): the two multipliers $\lambda$ and $\boldsymbol{\mu}$ measure how the observation $y_{obs}$ is good to reconstruct $y$ and $\pvec$, i.e. to satisfy the two constraints which appears in $\mathcal{V}$: $\rho_1^{-1}\mathcal{J}(y,\pvec)=0$ in $\Lvec^2(Q_T)$ and $\rho^{-1} \mathcal{I}(y,\pvec)=0$ in $L^2(Q_T)$.
\end{remark}

\

Eventually,  let us emphasize that, as in Section \ref{recovering_y} - the additional optimality system \eqref{eq:system_lambda_mu} can be used to define a equivalent saddle-point formulation.  Precisely, in view of \eqref{eq:system_lambda_mu}, we introduce the space $\Psi$ defined by 
\[
	\Psi=\{(\varphi,\boldsymbol{\sigma})\in 
	[C^0([0,T]; L^2(\Omega))\cap L^2(0,T;H^1_0(\Om))]\times\Lvec^2(Q_T),\,~\rho_0\mathcal{I}^\star(\varphi,\boldsymbol{\sigma})\in L^2(Q_T),
	\,~\varphi(\cdot,T)=0\}.
\]
	Endowed with the scalar product  
$$
\langle (\varphi,\boldsymbol{\sigma}),(\overline{\varphi},\overline{\boldsymbol{\sigma}})\rangle_{\Psi} := 
	\jjntQT\! \biggl(\boldsymbol{\sigma}\cdot\overline{\boldsymbol{\sigma}}+\rho^{-2}\nabla\varphi\cdot\nabla\overline{\varphi}+\rho^{2}_0\,\mathcal{I}^\star (\varphi,\boldsymbol{\sigma})
	\mathcal{I}^\star (\overline{\varphi},\overline{\boldsymbol{\sigma}})\biggr)\, dxdt,
$$
	we check that %, in view of assumption \eqref{iobs},
	  $\Psi$ is a Hilbert space.  Then, for any parameters $\mat{\alpha}=(\alpha_1,\alpha_2)\in (0,1)^2$
	  and $\mat{r}=(r_1,r_2)\in (\mathbb{R}^+_*)^2$, we consider the following mixed
	   formulation\,: find $((y,\pvec),(\varphi,\boldsymbol{\sigma}))\in \mathcal{U}\times \Psi$ such that 
\begin{equation} \label{eq:mfalphayp}
	\left\{
		\begin{array}{rcl}
			\noalign{\smallskip}\dis
			 a_{\boldsymbol{r},\mat{\alphasm}}((y,\pvec),(\ybar,\overline{\pvec})) + b_{\mat{\alphasm}}((\ybar,\overline{\pvec}),(\varphi,\boldsymbol{\sigma})) 
			 & = & l_{1,\mat{\alphasm}}(\ybar,\overline{\pvec}) 
			  \qquad \forall \ybar,\overline{\pvec} \in \mathcal{U}\\
			 \noalign{\smallskip}\dis
			 b_{\mat{\alphasm}}((y,\pvec), (\overline{\varphi},\overline{\boldsymbol{\sigma}}))-c_{\mat{\alphasm}}((\varphi,\boldsymbol{\sigma}),(\overline{\varphi},\overline{\boldsymbol{\sigma}}))
			 &=&l_{2,\mat{\alphasm}}(\overline{\varphi},\overline{\boldsymbol{\sigma}})
			  \qquad \forall \overline{\varphi},\overline{\boldsymbol{\sigma}}\in  \Psi,
		\end{array}
	\right.
\end{equation}
	where
\begin{align}
	& a_{\boldsymbol{r},\mat{\alphasm}} : \mathcal{U} \times\mathcal{U}\to \mathbb{R},\quad a_{\boldsymbol{r},\mat{\alphasm}}((y,\pvec),(\ybar,\overline{\pvec}))\!:=(1-\alpha_1)\!\jjntqT\!\!\!\rho_0^{-2}y\overline{y}\,dxdt
	+r_1\!\jjntQT\!\!\!\rho^{-2}_1\mathcal{J} ( y,\pvec)\cdot\mathcal{J} (\ybar,\overline{\pvec})\,dxdt\nonumber\\
	\noalign{\smallskip}\dis
	&\hspace{6.5cm}+r_2\jjntQT\rho^{-2}\mathcal{I}(y,\pvec)\mathcal{I}(\ybar,\overline{\pvec})\,dxdt,\nonumber\\
	\noalign{\smallskip}\dis
	& b_{\mat{\alphasm}}: \mathcal{U}\times \Psi \to \mathbb{R}, \quad b_{\mat{\alphasm}}((y,\pvec),(\varphi,\boldsymbol{\sigma})) := \jjntQT \mathcal{J}(y,\pvec)\cdot \boldsymbol{\sigma}\,dx\,dt+\jjntQT\mathcal{I}(y,\pvec) \varphi\,dx\,dt \nonumber\\
	\noalign{\smallskip}\dis
	&\hspace{6cm}-\alpha_1 \jjntqT  \mathcal{I}^\star (\varphi,\boldsymbol{\sigma})\,y\,dxdt, 
	\nonumber\\ 
	\noalign{\smallskip}\dis
	& c_{\mat{\alphasm}}: \Psi \times \Psi  \to \mathbb{R},  
	\quad c_{\mat{\alphasm}}((\varphi,\boldsymbol{\sigma}),(\overline{\varphi},\overline{\boldsymbol{\sigma}})) := 
	\alpha_1\jjntQT  \rho_0^2\mathcal{I}^\star (\varphi,\boldsymbol{\sigma})\mathcal{I}^\star 
	(\overline{\varphi},\overline{\boldsymbol{\sigma}})\, dxdt \nonumber\\
	\noalign{\smallskip}\dis
	&\hspace{6cm}+\alpha_2\jjntQT  \mathcal{J} (\varphi,\boldsymbol{\sigma})\cdot\mathcal{J} (\overline{\varphi},\overline{\boldsymbol{\sigma}})\, dxdt\nonumber\\
	\noalign{\smallskip}\dis
	& l_{1,\mat{\alphasm}}: \mathcal{U}\to \mathbb{R},  
	\quad
	l_{1,\mat{\alphasm}}(y,\pvec) := (1-\alpha_1)\jjntqT \rho^{-2}_0y_{obs} \, y\, dx dt, 
	\nonumber\\
	\noalign{\smallskip}\dis
	& l_{2,\mat{\alphasm}}: \Psi\to \mathbb{R},  
	\quad 
	l_{2,\mat{\alphasm}}(\varphi,\boldsymbol{\sigma}) := -\alpha_1 \jjntqT  y_{obs}\, \mathcal{I}^\star (\varphi,\boldsymbol{\sigma})\,dxdt.  \nonumber
\end{align}

	From the symmetry of $a_{\bf{r},{\mat{\mathsmaller{\alpha}}}}$ and $c_{{\mat{\mathsmaller{\alpha}}}}$, 
	we easily check that this formulation corresponds to the saddle point problem\,: 
\[
\left\{
	\begin{aligned}
	& \sup_{(\varphi,\boldsymbol{\sigma})\in \Psi} \inf_{(y,\pvec)\in \mathcal{U}} 
	\mathcal{L}_{\boldsymbol{r},\mat{\alphasm}}((y,\pvec),(\varphi,\boldsymbol{\sigma})), \\
	\noalign{\smallskip}\dis
	& \mathcal{L}_{\bf{r},\mat{\alphasm}}((y,\pvec),(\varphi,\boldsymbol{\sigma})):=\mathcal{L}_{\bf{r}}((y,\pvec),(\rho\varphi,\rho_1\boldsymbol{\sigma}))
	-\frac{\alpha_2}{2} \Vert \mathcal{J}(\varphi,\boldsymbol{\sigma})\Vert_{L^2(Q_T)}^2\\
	\noalign{\smallskip}\dis
	&\hspace{5.9cm}-\frac{\alpha_1}{2} \Vert \rho_0\mathcal{I}^\star(\varphi,\boldsymbol{\sigma}) +\rho_0^{-1}(y-y_{obs})1_{\om}\Vert_{L^2(Q_T)}^2. 
	\end{aligned}
\right.
\]
\begin{proposition}\label{prop:mfalphayp}
	Let $\rho_0\in \mathcal{R}$ and $\rho,\rho_1\in \mathcal{R}\cap L^{\infty}(Q_T)$. Then, for any $\alpha_1,\alpha_2\in (0,1)$, the formulation \eqref{eq:mfalphayp} is well-posed. Moreover, the unique pair 
	$((y,\pvec),(\varphi,\boldsymbol{\sigma}))$ in $\mathcal{U}\times \Psi$ satisfies 
\begin{equation}\label{estimate_solalpha_yp}
	\theta_1 \Vert(y,\pvec)\Vert_\mathcal{U}^2 + \theta_2 \Vert (\varphi,\boldsymbol{\sigma})\Vert_{\Psi}^2 
	\leq  \biggl(\frac{(1-\alpha_1)^2}{\theta_1} + \frac{\alpha_1^2}{\theta_2}\biggr)\Vert \rho_0^{-1}y_{obs}\Vert^2_{L^2(q_T)}. 
\end{equation}
	with 
\begin{equation}
	\theta_1:=\min\biggl(1-\alpha_1, \frac{r_1}{\eta_1}, \frac{r_2}{\eta_2}\biggr) , \quad \theta_2:=C(\mat{\alpha},\rho_\star).   \nonumber
\end{equation}
\end{proposition}
\textsc{Proof-} We easily get the continuity of the bilinear forms $a_{\bf{r},\mat{\alphasm}}$, $b_{\mat{\alphasm}}$ and $c_{\mat{\alphasm}}$: 
 $$
	\begin{aligned}
		& \vert a_{\boldsymbol{r},\mat{\alphasm}}((y,\pvec),(\ybar,\overline{\pvec}))\vert 
		\leq \max\{1-\alpha_1,r_1\eta_1^{-1},r_2\eta_2^{-1}\} \Vert (y,\pvec)\Vert_\mathcal{U} 
		\Vert(\ybar,\overline{\pvec})\Vert_\mathcal{U}, \\
		\noalign{\smallskip}\dis
		& \vert b_{\mat{\alphasm}}((y,\pvec), (\varphi,\boldsymbol{\sigma}))\vert  \leq \max\{\alpha_1,\eta_1^{-1/2}\|\rho_1\|_{L^{\infty}
		(Q_T)},\eta_2^{-1/2}\|\rho\|_{L^{\infty}(Q_T)}\} 
		\Vert (y,\pvec)\Vert_\mathcal{U} \Vert (\varphi,\boldsymbol{\sigma})\Vert_{\Psi} ,\\
		\noalign{\smallskip}\dis
		& \vert c_{\mat{\alphasm}}((\varphi,\boldsymbol{\sigma}),(\overline{\varphi},\overline{\boldsymbol{\sigma}}))\vert 
		\leq \max\{\alpha_1,\alpha_2\} \Vert (\varphi,\boldsymbol{\sigma})\Vert_{\Psi} \Vert(\overline{\varphi},\overline{\boldsymbol{\sigma}})\Vert_{\Psi} 
	\end{aligned}
 $$
 	for all $(y,\pvec),(\ybar,\overline{\pvec})\in \mathcal{U}$ and 
	for all $(\varphi,\boldsymbol{\sigma}),(\overline{\varphi},\overline{\boldsymbol{\sigma}})\in \Psi$. Also, we can easily deduce 
	the continuity of the linear form $l_{1,\mat{\alphasm}}$ and $l_{2,\mat{\alphasm}}$ : 
	$\Vert l_{1,\mat{\alphasm}}\Vert_{\mathcal{U}^{\prime}}\le (1-\alpha_1)\Vert\rho_0^{-1} y_{obs}\Vert_{L^2(q_T)}$ and
	$\Vert l_{2,\mat{\alphasm}}\Vert_{\Psi^{\prime}} \le \alpha_1\Vert\rho_0^{-1} y_{obs}\Vert_{L^2(q_T)}$.   
 
 	Moreover, we also obtain the coercivity of $a_{\boldsymbol{r},\mat{\alphasm}}$ and of $c_{\mat{\alphasm}}$: precisely, we check 
	that $$a_{\boldsymbol{r},\mat{\alphasm}}((y,\pvec),(y,\pvec)) \geq \theta_1 \Vert (y,\pvec)\Vert^2_\mathcal{U}
	\quad\forall (y,\pvec)\in \mathcal{U}$$ 
	while, for any $m\in (0,1)$, denote $g=\rho_0\mathcal{I}^\star (\varphi,\boldsymbol{\sigma})$ 
	and $\Gvec= \mathcal{J} (\varphi,\boldsymbol{\sigma})$, by writing, 
\[
	\begin{aligned}
		c_{\mat{\alphasm}}((\varphi,\boldsymbol{\sigma}),(\varphi,\boldsymbol{\sigma})) = &~ 
		\alpha_1 \Vert g\Vert^2_{L^2(Q_T)}
		+\alpha_2 \Vert \Gvec\Vert^2_{L^2(Q_T)}\\
		=&~
		m(\alpha_1 \Vert g\Vert^2_{L^2(Q_T)}
		+\alpha_2 \Vert \Gvec\Vert^2_{L^2(Q_T)})+C_{\Om,T}(1-m)\min\{\alpha_1\rho_\star,\alpha_2\} \|\nabla\varphi\|_{L^2(Q_T)}^2\\
		\ge &~ 
		m\alpha_1\Vert g\Vert^2_{L^2(Q_T)}
		+\min\{m\alpha_2,{C_{\Om,T}(1-m)\min\{\alpha_1\rho_\star,\alpha_2\}\over2} \}\|\boldsymbol{\sigma}\|_{L^2(Q_T)}^2\\
		&+{(1-m)\over2}\min\{\alpha_1\rho_\star,\alpha_2\} \|\nabla\varphi\|_{L^2(Q_T)}^2
	\end{aligned}
\]
	we get $c_{\mat{\alphasm}} ((\varphi,\boldsymbol{\sigma}),(\varphi,\boldsymbol{\sigma})) \geq \theta_2\|(\varphi,\boldsymbol{\sigma})\|_{\Psi}^2$ for all 
	$(\varphi,\boldsymbol{\sigma})\in \Psi$.

The result \cite[Prop 4.3.1]{brezzi_new}  implies the well-posedness of the mixed formulation \eqref{eq:mfalphayp} and the estimate \eqref{estimate_solalpha_yp}. \Fin

	The $\mat{\alpha}$-term in $\mathcal{L}_{\bf{r},{\mat{\mathsmaller{\alpha}}}}$ is a stabilization term\,: it ensures 
	a coercivity property of $\mathcal{L}_{\bf{r},{\mat{\mathsmaller{\alpha}}}}$ with respect to the variables 
	$(\varphi,\sigma)$ and automatically the well-posedness. In particular, there is no need to prove any inf-sup 
	property for the bilinear form $b_{{\mat{\mathsmaller{\alpha}}}}$.

\begin{proposition}
	If the solution $((y,\pvec),(\lambda,\mu)\in \mathcal{U}\times L^2(Q_T)\times \Lvec^2(Q_T)$ 
	of \eqref{eq:mf4} enjoys the property $(\lambda,\mu)\in \Psi$, then the solutions of
	 \eqref{eq:mf4} and \eqref{eq:mfalphayp} coincide. 
\end{proposition}
\textsc{Proof-} The hypothesis of regularity and the relation \eqref{eq:system_lambda_mu} imply that the solution 
	$((y,\pvec),(\lambda,\mu)\in \mathcal{U}\times L^2(Q_T)\times \Lvec^2(Q_T)$ of \eqref{eq:mf4} is also a solution of \eqref{eq:mfalphayp}. 
	The result then follows from the uniqueness of the two formulations. \Fin

\subsection{Dual formulation of the extremal problem (\ref{eq:mf4})}

For any $\boldsymbol{r}=(r_1,r_2)\in (\mathbb{R}^+)^2$, we define the linear operator $\mathcal{T}_{\boldsymbol{r}}$ from $\mathcal{X}=L^2(Q_T)\times \Lvec^2(Q_T)$ into $\mathcal{X}$ by
\begin{equation}
\mathcal{T}_{\boldsymbol{r}}(\lambda,\boldsymbol{\mu}):=(\rho_1^{-1} \mathcal{J}(y,\pvec),\rho^{-1}\mathcal{I}(y,\pvec))\nonumber
\end{equation}
where $(y,\pvec)\in \mathcal{U}$ solves, for any $\boldsymbol{r}=(r_1,r_2)\in \mathbb{R}^*_+$  
\begin{equation}\label{eqTr}
a_{\boldsymbol{r}}((y,\pvec),(\ybar,\overline{\pvec}))=b((\ybar,\overline{\pvec}),(\lambda,\boldsymbol{\mu})) \quad  \forall\,(\ybar,\overline{\pvec})\in \mathcal{U}.
\end{equation}

Similarly to Lemma \ref{lemmaA}, the following holds true. 

\begin{lemma}
For any $\boldsymbol{r}=(r_1,r_2)\in \mathbb{R}^*_+$, the operator $\mathcal{T}_{\boldsymbol{r}}$ is a strongly elliptic, symmetric isomorphism from $\mathcal{X}$ into $\mathcal{X}$.
\end{lemma}
\textsc{Proof-} From the definition of $a_{\boldsymbol{r}}$, we get that $\Vert \mathcal{T}_{\boldsymbol{r}}(\lambda,\boldsymbol{\mu})\Vert_{\mathcal{X}}\leq \min(r_1,r_2)^{-1}\Vert(\lambda,\boldsymbol{\mu})\Vert_{\mathcal{X}}$ leading to the continuity of $\mathcal{T}_{\boldsymbol{r}}$. Next, consider any $(\lambda^{\prime},\boldsymbol{\mu}^{\prime})\in \mathcal{X}$ and denote by $(y^{\prime},\pvec^{\prime})$ the corresponding solution of (\ref{eqTr}) so that $\mathcal{T}_{\boldsymbol{r}}(\lambda^{\prime},\boldsymbol{\mu}^{\prime})=(\rho_1^{-1} \mathcal{J}(y^{\prime},\pvec^{\prime}),\rho^{-1}\mathcal{I}(y^{\prime},\pvec^{\prime}))$. Relation (\ref{eqTr}) with $(\ybar,\overline{\pvec})=(y^{\prime},\pvec^{\prime})$ implies that 
\begin{equation}
\jjntQT \mathcal{T}_{\boldsymbol{r}}(\lambda^\prime,\boldsymbol{\mu}^{\prime})\cdot (\lambda,\boldsymbol{\mu})\,dxdt=a_{\boldsymbol{r}}((y,\pvec),(y^{\prime},\pvec^{\prime})) \label{arAlambdabis}
\end{equation}
and therefore the symmetry and positivity of $\mathcal{T}_{\boldsymbol{r}}$. The last relation with $(\lambda^{\prime},\boldsymbol{\mu}^{\prime})$ together with the observability estimate (\ref{crucial_estimate_good_fo}) imply that the operator $\mathcal{T}_{\boldsymbol{r}}$ is also positive definite. 
Actually, as announced, we can check that $\mathcal{T}_{\boldsymbol{r}}$ is strongly elliptic, i.e. there exits a constant $C>0$ such that 
\begin{equation} \nonumber
\jjntQT  \mathcal{T}_{\boldsymbol{r}}(\lambda,\boldsymbol{\mu})\cdot (\lambda,\boldsymbol{\mu})  \, dxdt\geq C \Vert (\lambda,\boldsymbol{\mu})\Vert^2_{\mathcal{X}}, \quad\forall (\lambda,\boldsymbol{\mu})\in \mathcal{X}.
\end{equation}
We argue by contradiction and suppose that there exists a sequence $\{(\lambda_n,\boldsymbol{\mu_n})\}_{n\geq 0}$ of $\mathcal{X}$ such that 
\begin{equation}\label{hyp_contr}
\Vert (\lambda_n,\boldsymbol{\mu_n})\Vert_{\mathcal{X}}=1 \quad \forall n\geq 0 \quad \textrm{and}\quad\lim_{n\to \infty} \jjntQT  \mathcal{T}_{\boldsymbol{r}}(\lambda_n,\boldsymbol{\mu_n})\cdot (\lambda_n,\boldsymbol{\mu_n})\, dxdt =0.
\end{equation}
We denote by  $(y_n,\pvec_n)$ the solution of (\ref{eqTr}) corresponding to $(\lambda_n,\boldsymbol{\mu_n})$. From (\ref{arAlambdabis}), we then obtain that 
\begin{equation}
\lim_{n\to} \biggl(\Vert \rho_0^{-1} y_n\Vert^2_{L^2(q_T)} + r_1 \Vert\rho_1^{-1} \mathcal{J}(y_n,\pvec_n)\Vert^2_{\Lvec^2(Q_T)}+r_2 \Vert\rho^{-1} \mathcal{I}(y_n,\pvec_n)\Vert^2_{L^2(Q_T)}\biggr)=0.  \label{limitbis}
\end{equation}
Moreover, from (\ref{eqTr}) with $(y,\pvec)=(y_n,\pvec_n)$ and $(\lambda,\boldsymbol{\mu})=(\lambda_n,\boldsymbol{\mu_n})$, we get the equality
\begin{equation}
\label{26bis}
\begin{aligned}
&\jjntQT \biggl(\rho_1^{-1} \mathcal{J}(\ybar,\overline{\pvec})\cdot (r_1 \rho_1^{-1}\mathcal{J}(y_n,\pvec_n)-\boldsymbol{\mu_n}) + \rho^{-1} \mathcal{I}(\ybar,\overline{\pvec})\dot (r_2 \rho^{-1}\mathcal{I}(y_n,\pvec_n)-\lambda_n)\biggr)\,dxdt \\
& +\jjntqT \rho_0^{-2}y_n\ybar \,dxdt=0
\end{aligned}
\end{equation}
for every $(\ybar,\overline{\pvec}) \in \mathcal{U}$. Then, in order to get a contradiction, we define the sequence $\{(\ybar_n,\overline{\pvec_n})\}_{n\geq 0}$ as follow: 

\begin{equation}
\nonumber%\label{eq:heat_mixed}
	\left\{
		\begin{array}{lll}
		      \rho_1^{-1} \mathcal{J}(\ybar_n,\overline{\pvec_n})=r_1 \rho_1^{-1}\mathcal{J}(y_n,\pvec_n)-\boldsymbol{\mu_n} & \textrm{in} & Q_T \\
		      \rho^{-1} \mathcal{I}(\ybar_n,\overline{\pvec_n})=r_2 \rho^{-1}\mathcal{I}(y_n,\pvec_n)-\lambda_n & \textrm{in} & Q_T \\
			\ybar_n  = 0										  		&\textrm{on}& \Sigma_T, 	\\
   			\ybar_n(x, 0) = 0 								  		& \textrm{in}& \Om,
   		\end{array} 
 	\right.
\end{equation}
so that, for all n, $(\ybar_n,\overline{\pvec_n})$ is the solution of a first order system as discussed in the Appendix with zero initial data and source term in $\mathcal{X}$. Energy estimates (\ref{energy_est_mixed}) implies that 
\begin{equation}
\nonumber
\begin{aligned}
\Vert \rho_0^{-1}\ybar_n\Vert_{L^2(q_T)} \leq C_{\Omega,T} \rho_\star^{-1} \biggl(  \Vert \rho_1&\Vert_{L^{\infty}(Q_T)} \Vert r_1 \rho_1^{-1}\mathcal{J}(y_n,\pvec_n)-\boldsymbol{\mu_n} \Vert_{\Lvec^2(Q_T)} \\
& + \Vert \rho\Vert_{L^{\infty}(Q_T)} \Vert r_2 \rho^{-1}\mathcal{I}(y_n,\pvec_n)-\lambda_n\Vert_{L^2(Q_T)} \biggr)
\end{aligned}
\end{equation}
for some constant $C_{\Omega,T}$
and that $(\ybar_n,\overline{\pvec_n})\in \mathcal{U}$. Then, using this inequality and (\ref{26bis}) with $(\ybar,\overline{\pvec})=(\ybar,\overline{\pvec}_n)$, we get that 
\begin{equation}
\nonumber
\begin{aligned}
\Vert r_1 \rho_1^{-1}\mathcal{J}(y_n,\pvec_n)-\boldsymbol{\mu_n} &\Vert_{\Lvec^2(Q_T)}+ \Vert r_2 \rho^{-1}\mathcal{I}(y_n,\pvec_n)-\lambda_n \Vert_{L^2(Q_T)} \\
&\leq 2 C_{\Omega,T} \rho_\star^{-1} \max\biggl(\Vert \rho_1\Vert_{L^{\infty}(Q_T)},\Vert \rho\Vert_{L^{\infty}(Q_T)}\biggr) \Vert\rho_0^{-1} y_n\Vert_{L^2(q_T)}.
\end{aligned}
\end{equation}
Eventually, from (\ref{limitbis}), we conclude that $\lim_{n\to \infty} \Vert \lambda_n\Vert_{L^2(Q_T)}=\lim_{n\to \infty} \Vert \boldsymbol{\mu_n}\Vert_{\Lvec^2(Q_T)}=0$, which is in contradiction with the first hypothesis of (\ref{hyp_contr}).
\Fin

Again, the introduction of the operator $\mathcal{T}_{\boldsymbol{r}}$ is motivated by the following proposition, which reduces the determination of the solution $(y,\pvec)$ of Problem (\ref{inverse_pb_2}) to the unconstrained minimization of a elliptic functional. 

\begin{proposition}\label{prop_equiv_dual_mixed}
For any $\boldsymbol{r}>0$, let $(y_0,\pvec_0)\in \mathcal{U}$ be the unique solution of 
\[
a_{\boldsymbol{r}}((y_0,\pvec_0),(\ybar,\overline{\pvec})= l(\ybar,\overline{\pvec}), \quad \forall (\ybar,\overline{\pvec})\in \mathcal{U}
\]
and let $J_{\boldsymbol{r}}^{\star\star}:\mathcal{X}\to \mathcal{X}$ be the functional defined by 
\[
J_{\boldsymbol{r}}^{\star\star}(\lambda,\boldsymbol{\mu}) = \frac{1}{2} \jjntQT \mathcal{T}_r(\lambda,\boldsymbol{\mu})\cdot (\lambda,\boldsymbol{\mu}) \,dx\,dt - b((y_0,\pvec_0), (\lambda,\boldsymbol{\mu})).
\]
The following equality holds : 
\begin{equation}\nonumber
\sup_{(\lambda,\boldsymbol{\mu})\in \mathcal{X}}\inf_{(y,\pvec)\in \mathcal{U}} \mathcal{L}_{\boldsymbol{r}}((y,\pvec),(\lambda,\boldsymbol{\mu})) = - \inf_{(\lambda,\boldsymbol{\mu})\in \mathcal{X}} J_r^{\star\star}(\lambda,\boldsymbol{\mu})\quad + \mathcal{L}_r((y_0,\pvec_0),(0,\boldsymbol{0})).
\end{equation}
where the Lagrangien $\mathcal{L}_{\boldsymbol{r}}$ is defined in Remark \ref{rk6}.
\end{proposition}

\section{Concluding remarks and perspectives}  \label{sec_conclusion}

The mixed formulations we have introduced in order to address inverse problems for linear parabolic type equations seems original. These formulations are nothing else than the Euler systems associated to weighted least-squares type functionals and depend on both the state to reconstruct and a Lagrange multiplier. This multiplier is introduced to take into account the state constraint $Ly-f=0$ and turns out to be a measure of how good the observation data is to reconstruct the solution.  This approach, recently used in a controllability context in \cite{munch_desouza}, leads to a variational problems defined over time-space functional Hilbert spaces, without distinction between the time and the space variable. 
The main ingredient is the unique continuation property leading to well-posedness in appropriate constructed Hilbert spaces. Moreover, global Carleman estimates then allow to precise in which norm the full solution can be reconstructed.  For these reasons, the method can be applied to many systems for which such estimates are available, as in \cite{NC-AM-InverseProblems} for linear hyperbolic equations,  or as in Section 
\ref{recovering_y_fo} for a first order system. In the parabolic situation, in view of regularization property, the method requires the introduction of exponentially vanishing weights at the initial time: this guarantees a stable Lipschitz reconstruction of the solution on the whole domain, the initial condition excepted. 

On the theoretical standpoint, the minimization of the $L^2$-weighted least-squares norm with respect either to $y\in \mathcal{W}$ (Problem \ref{P}, Section \ref{sec1_direct}), either to the initial data $y_0\in \mathcal{H}$ (Problem \ref{extremal_problem}, Section \ref{sec:intro}) is equivalent. However, the completed space $\mathcal{W}$ embedded in the space $C([\delta,T],H^0_1(\Omega))$ is a priori much more 'practical' than the huge space $\mathcal{H}$, a fortiori since from the definition of the cost, the variable of interest is not $y$ but $\rho_0^{-1}y\in C([0,T],H_0^1(\Omega))$ with $\rho_0^{-1}(\cdot,t=0)=0$ in $\Omega$. Therefore, on a practical (i.e. numerical) viewpoint, as enhanced in \cite{imanuv_puel_yama_2,puel_tycho} and recently used in \cite{bourgeois_parabolic,BeilinaKlibanov14} for inverse problems  and in \cite{EJAM_pedregal,munch_desouza} in the close controllability context (see Remark \ref{rk_control_analogy}), variational methods where the state $y$ is kept as the unknown are very appropriate and lead to robust approximations. Moreover, as detailed in \cite{NC-AM-InverseProblems}, the space-time framework allows to use classical approximation and interpolation theory leading to strong convergence results with error estimates, again without the need of proving any discrete Carleman inequalities. 
We refer to the second part \cite{munch-desouza-Part2} of this work  where the numerical approximation of the mixed formulations (\ref{eq:mf}) in $(y,\lambda)$ and (\ref{eq:mf4}) in $((y,\pvec),(\lambda,\boldsymbol{\mu}))$ is examined, implemented and compared with the standard minimization of the cost with respect to the initial data. As observed in \cite{munch_desouza} section 3.2 for the related control problem, described in Remark \ref{rk_control_analogy}, an appropriate preliminary change (renormalization) of variable, i.e. $\tilde{y}:=\rho_0^{-1}y$, so as to eliminate (by compensation) the exponential behavior of the coefficient in $\rho^{-1}Ly=\rho^{-1}L(\rho_0 \tilde{y})$, leads to an impressive low condition number of the corresponding discrete system. We also emphasize, that the second mixed formulation (\ref{eq:mf4}), apparently more involved with more variables allows to use (standard) continuous finite dimensional approximation spaces for $\mathcal{U}$, in contrast to the formulation (\ref{eq:mf})
which requires continuously differentiable approximation spaces.

Eventually, we also emphasize that such direct method may be used to reconstruct the state as well as a source term. Assuming that the source $f(x,t)=\sigma(t) \mu(x)$ with $\sigma\in C^1([0,T])$, $\sigma(0)\neq 0$ and $\mu\in L^2(\Omega)$, it is shown in \cite{choulli_yamamoto} that the knowledge of $\partial_t(\partial_{\nu}u)\in L^2(\partial\Omega\times (0,T))$ allows to reconstruct uniquely the pair $(y,\mu)$ satisfying the state equation $Ly-\sigma \mu=0$. This allows to construct appropriate Hilbert spaces, associate a least-squares functional in $(y,\mu)$ and the corresponding optimality system. The (logarithmic) stability estimate proved in \cite[Theorem 1.2]{choulli_yamamoto} guarantees the reconstruction of the solution. 
We refer to the Section 3 of \cite{NC-AM-InverseProblemsBoundary} where this strategy is implemented in the simpler case of the wave equation.

\begin{appendix} \label{parabolic-mixed-form}

\section[Well-posedness of Parabolic equations in the mixed form]{Appendix: Well-posedness of Parabolic equations in the mixed form}

The aim of this appendix is to study the existence and uniqueness of weak solution and solution by transposition for the following linear boundary value problem, which appears in Section \ref{recovering_y_fo} : find $(y,\pvec)$ such that 
\begin{equation}
\label{eq:heat_mixed_form}
	\left\{
		\begin{array}{lll}
   			y_t- \nabla\cdot\pvec+d\,y  =f, \quad c(x)\nabla y-\pvec=\Fvec 	& \textrm{in}& Q_T,    		\\
   			y  = 0										  		&\textrm{on}& \Sigma_T, 	\\
   			y(x, 0) = y_0(x) 								  		& \textrm{in}& \Om.
   		\end{array} 
 	\right.
\end{equation}
We assume that the initial datum $y_0$ belongs to $L^2(\Om)$ and that the source terms $f$ and $\Fvec$ belong to $L^2(Q_T)$ and 
	$\Lvec^2(Q_T)$, respectively. The functions $c$ and $d$ enjoys the regularity described in the introduction: $c:=(c_{i,j})\in C^1(\overline\Om;\mathcal{M}_{N}(\mathbb{R}))$ with $(c(x)\xi,\xi)\geq c_0|\xi|^2$ for any $x\in\overline\Om$ and $\xi\in \mathbb{R}^N$~$(c_0>0)$ and $d \in L^\infty(Q_T)$.
		
	First, let us introduce a definition of weak solution in accordance to the classical definition 
	of weak solution for the \textit{standard} parabolic equation (\ref{eq:heat}).
\begin{definition}\label{def:weak_heat_stand_mixed}
	We say that a pair $(y,\pvec)$, satisfying
\begin{equation}	
\label{eq:reg_heat_mixed}
   	\pvec\in \Lvec^2(Q_T),\,~ y\in L^2(0,T;H^1_0(\Omega)),\,~\hbox{with}\,~y_t\in L^2(0,T;H^{-1}(\Om)),
\end{equation}
	is a weak solution of the parabolic equation in the mixed form \eqref{eq:heat_mixed_form} if and only if\,:
\begin{itemize}
	\item [$(i)$] $\langle y_t,w\rangle_{H^{-1}(\Om),H^1_0(\Om)}+(\pvec,\nabla w) 
	+ (d\,y,w)=(f,w)$ for all $w\in H^1_0(\Om)$ and a.e. time $t\in[0,T]$;
	\item [$(ii)$] $ (\nabla y,\uvec)-(c^{-1}\pvec, \uvec) =(c^{-1}\Fvec,\uvec)$ for all $\uvec\in \Lvec^2(\Om)$ and a.e. time $ t\in[0,T]$;
\end{itemize}
\begin{enumerate}
	\item [(iii)] $y(\cdot,0)=y_0$.
\end{enumerate}
Here, $(\cdot,\cdot)$ denotes the inner product in $L^2(\Omega)$ and in $\Lvec^2(\Omega)$. 
\end{definition}

\begin{remark}\label{rmq_mixed}
According to 	Definition \eqref{def:weak_heat_stand_mixed}, a weak solution for (\ref{eq:heat_mixed_form}) is a weak solution for the standard
	parabolic equation\,:
\begin{equation}
\label{eq:parab}
	\left\{
		\begin{array}{lll}
   			y_t- \nabla\cdot(c(x)\nabla y)+d\,y  =f-\nabla\cdot\Fvec 		& \textrm{in}& Q_T,    		\\
   			y  = 0										  		&\textrm{on}& \Sigma_T, 	\\
   			y(x, 0) = y_0(x) 								  		& \textrm{in}& \Om.
   		\end{array} 
 	\right.
\end{equation}
	Indeed, taking $\uvec=c\,\nabla w$ and summing the equations in $(i)$ and $(ii)$, we obtain the definition of weak solution for \eqref{eq:parab}.
\end{remark}

\begin{proposition}\label{ex_un_mixed_sol1}
	There exists a unique weak solution for the parabolic equation in the mixed form~\eqref{eq:heat_mixed_form}. 
	Moreover, there exists a constant $C>0$ such that
\begin{equation}\label{energy_est_mixed}
	\begin{alignedat}{2}
	\|y'\|_{L^2(0,T;H^{-1}(\Omega))}+\|y\|_{L^2(0,T;H^1_0(\Omega))}+\|\pvec\|_{\Lvec^2(Q_T)}\leq C(\|y_0\|_{L^2(\Omega)}+\|f\|_{L^2(Q_T)}+\|\Fvec\|_{\Lvec^2(Q_T)}).
	\end{alignedat}
\end{equation}
\end{proposition}
\textsc{Proof-} Following \cite{jerome_space_time}, the proof of existence of solution relies on the {\it Faedo-Galerkin method} and it is divided in several steps.
	\begin{itemize}
	\item[a)] {\bf Galerkin approximations.} We first introduce some notations : let $\{w_k\,:\,k\in\mathbb{N}\}$ be a
	orthogonal basis of $H^1_0(\Omega)$ (which is orthonormal in $L^2(\Omega)$) and 
	$\{\uvec_k\,:\,k\in\mathbb{N}\}$ be a orthonormal basis of $\Lvec^2(\Omega)$. For each pair 
	$(m,n)\in \mathbb{N}\times \mathbb{N}$, we look for a pair 
	$(y_n,\pvec_m):[0,T]\to H^1_0(\Omega)\times \Lvec^2(\Omega)$ of the form
\begin{equation}\label{eq:y_p_n}
	y_n(t)=\sum_{k=1}^na_n^k(t)w_k\quad\hbox{and}\quad \pvec_m(t)=\sum_{k=1}^{m}b_m^k(t)\uvec_k,
\end{equation}
solution of the weak formulation : 
\begin{equation}\label{eq:y_and_p}
	\left\{
		\begin{array}{ll}
			<y^\prime_n,w_i>+(\pvec_m,\nabla w_i) 
			+ (d\,y_n,w_i)=(f,w_i)\quad (0\leq t\leq T,\,\,i=1,\ldots,n),\\
			 \noalign{\smallskip}\dis
			 (\nabla y_n,\uvec_j)-(c^{-1}\pvec_m,\uvec_j) =(c^{-1}\Fvec,\uvec_j)\quad (0\leq t\leq T,\,\,j=1,\ldots,m),
		\end{array}
		\right.
\end{equation}	
(the prime $^\prime$ stands for the derivation in time). $a_n^k$ ($k=1,\dots,n$) and $b_m^k$ ($k=1,\dots,m$) denote some time functions from $[0,T]$ to $\mathbb{R}$ for each pair $(m,n)$. We assume that the $a_n^k$ satisfies
\begin{equation}\label{eq:y_p_0}
	a_n^k(0)=(y_0,w_k)\quad (k=1,\ldots,n).
\end{equation}
We also note 
\[
	\begin{aligned}
			& (\fvec_n(t))_i=(f(t),w_i),\,~(\Yvec_0)_i=(y_0,w_i),\\
			& (\Avec_n)_{ij}=(w_j,w_i),\,~
			(\Dvec_n(t))_{ij}=(d(\cdot,t)\,w_j,w_i),
	\end{aligned}
	\qquad \forall i,j=1,\ldots,n
\]
and
\[
	\begin{array}{ll}
			 (\Bvec_m)_{ij}=(c^{-1}\uvec_j,\uvec_i),\quad(\Fvec_m(t))_i=(c^{-1}\Fvec(t),\uvec_i),\quad\forall i,j=1,\ldots,m
	\end{array}
\]
	and
\[
	\begin{array}{ll}
			  (\Evec_{nm})_{ij}=(\uvec_j,\nabla w_i) 
			  \quad \forall i=1,\ldots,n\quad \hbox{and}\quad j=1,\ldots,m.
	\end{array}
\]
Eventually, we also denote by $\Yvec_n(t)$ the vector formed by $a_n^k$ ($k=1,\ldots,n$) and
	$\Pvec_m(t)$ the vector formed by $b_m^k$  ($k=1,\ldots,m$).
	With thes notations, \eqref{eq:y_and_p} may be rewritten as
\begin{equation}\label{eq:y_and_p_first}
	\left\{
		\begin{array}{ll}
			\Yvec_n'(t)+\Evec_{nm}\Pvec_m(t)+\Dvec_n(t)\Yvec_n(t)=\fvec_n(t),\\
			 \noalign{\smallskip}\dis
			\Evec_{nm}^T\Yvec_n(t)-\Bvec_m\Pvec_m(t) =\Fvec_m(t),\\
			 \noalign{\smallskip}\dis
			 \Yvec_n(0)=\Yvec_0.
		\end{array}
		\right.
\end{equation}

From the positivity of $c$, $\Bvec_m$ is a symmetric and positive definite square matrix of order $m$: therefore $\Bvec_m$ is invertible and the second equation of \eqref{eq:y_and_p_first} implies the relation  
$$
\Pvec_m(t) =\Bvec_m^{-1}(\Evec_{nm}^T\Yvec_n(t)-\Fvec_m(t)).
$$
Using this equality in the first equation of \eqref{eq:y_and_p_first}, we obtain
\begin{equation}\label{eq:y_and_p_second}
	\left\{
		\begin{array}{ll}
			\Yvec_n'(t)+(\Evec_{nm}\Bvec_m^{-1}\Evec_{nm}^T+\Dvec_n(t))\Yvec_n(t)=\fvec_n(t)+
			\Evec_{nm}\Bvec_m^{-1}\Fvec_m(t),  \quad \textrm{for a.e. }t\in [0,T]\\
			 \noalign{\smallskip}\dis
			 \Yvec_n(0)=\Yvec_0.
		\end{array}
		\right.
\end{equation}

	\eqref{eq:y_and_p_second} is a system of $n$ linear ODEs of order $1$: hence, from standard theory for ODEs, there exists a unique absolutely continuous $\Yvec_n:[0,T]\to \mathbb{R}^n$ 
	satisfying \eqref{eq:y_and_p_second}. Consequently, the pair $(y_n,\pvec_m)$ given by \eqref{eq:y_p_n} is the unique solution of \eqref{eq:y_and_p} and \eqref{eq:y_p_0}.
	
	\item[b)] {\bf A priori estimates.} We now derive some uniform estimates for the pair $(y_n, \pvec_m)$
	with respect to $m$ and $n$. This allows to see that the sequences $\{y_n\}_{n>0}$, $\{\pvec_m\}_{m>0}$ converge to $y$ and $\pvec$ respectively, such that $(y,\pvec)$ is the weak solution of \eqref{eq:heat_mixed_form}.

	Multiplying the first equation of \eqref{eq:y_and_p} by $a_n^k(t)$ and summing over $k=1,\ldots, n$ 
	and the second equation of \eqref{eq:y_and_p} by $-b_m^k(t)$ and summing over $k=1,\ldots, m$, we obtain the relations
\[
	\left\{
		\begin{array}{ll}
			<y'_n,y_n>+(\pvec_m,\nabla y_n) 
			+ (d\,y_n,y_n)=(f,y_n)\quad (0\leq t\leq T),\\
			 \noalign{\smallskip}\dis
			 -(\nabla y_n,\pvec_m)+(c^{-1}\pvec_m,\pvec_m) =-(c^{-1}\Fvec,\pvec_m)\quad (0\leq t\leq T).
		\end{array}
		\right.
\]

Adding these two equations and applying the Cauchy-Schwarz inequality, we get\,:
\[
\begin{aligned}
	& \dis{d\over dt}\| y_n(t)\|_{L^2(\Omega)}^2+\|\pvec_m(t)\|_{\Lvec^2(\Omega)}^2\\
	&\leq C\biggl(\|y_n(t)\|_{L^2(\Omega)}^2+\|f(\cdot,t)\|_{L^2(\Omega)}\|y_n(t)\|_{L^2(\Omega)}+\|\Fvec(\cdot,t)\|_{\Lvec^2(\Omega)}\|\pvec_m(t)\|_{\Lvec^2(\Omega)}\biggr)
	\end{aligned}
\]
for some positive constant $C=C(\Vert c\Vert_{C^1(\overline{\Omega})},\Vert d\Vert_{L^{\infty}(Q_T)})$.
	Then, by using the Gronwall's Lemma, we deduce, that for all $m,n>0$
\begin{equation}\label{est_y_infty}
	\| y_n\|^2_{L^\infty(0,T;L^2(\Om))}\leq C(\|y_0\|_{L^2(\Omega)}^2+\|f\|_{L^2(Q_T)}^2+\|\Fvec\|_{\Lvec^2(Q_T)}^2)
\end{equation}
	and 
\begin{equation}\label{est_p_l2}
	\| \pvec_m\|^2_{\Lvec^2(Q_T)}\leq C(\|y_0\|_{L^2(\Omega)}^2+\|f\|_{L^2(Q_T)}^2+\|\Fvec\|_{\Lvec^2(Q_T)}^2).
\end{equation}

	Now we derive a uniform estimate for $\partial_t y_n$. To do this, fix any $w\in H^1_0(\Om)$, with 
	$\|w\|_{H^1_0(\Om)}\leq 1$. Notice that we can decompose $w$ as $w=w^1+w^2$ with $w^1\in span\{w_k\}_{k=1}^n$ 
	and $(w^2,w_k)=0$ for $k=1,\ldots,n$. Using
	the first equation of \eqref{eq:y_and_p}, we can deduce, for a.e. $0\leq t\leq T$, that
\[
			<y'_n,w^1>+(\pvec_m,\nabla w^1) 
			+ (d\,y_n,w^1)=(f,w^1).
\]
	
	Then, using that $(w^2,w_k)=0$ for $k=1,\ldots,n$ and $y_n'(t)=\sum_{k=1}^n(a_n^k)'(t)w_k$, we write
\[
			<y'_n,w>=( y'_n,w^1+w^2)=( y'_n,w^1)=(f,w^1)-(\pvec_m,\nabla w^1) - (d\,y_n,w^1).
\]
	Consequently,
\[
	\begin{aligned}
			|< y'_n(t),w>|\leq&~\|f(\cdot,t)\|_{L^2(\Omega)}\|w^1\|_{L^2(\Omega)}+\|\pvec_m(t)\|_{\Lvec^2(\Omega)}\|\nabla w^1\|_{\Lvec^2(\Omega)} \\
			&\hspace{1cm}+\|d\|_{L^\infty(Q_T)}\|y_n(t)\|_{L^2(\Omega)}\|w^1\|_{L^2(\Omega)}\\
			\leq&~C(\|f(\cdot,t)\|_{L^2(\Omega)}+\|\pvec_m\|_{\Lvec^2(\Omega)} +\|y_n\|_{L^2(\Omega)}),
	\end{aligned}
\]
	since $\|w^1\|_{H^1_0(\Omega)}\leq 1$. Finally, using \eqref{est_y_infty} and \eqref{est_p_l2}, we get
\begin{equation}\label{est_y_t}
	\begin{alignedat}{2}
	\|y'_n\|^2_{L^2(0,T;H^{-1}(\Omega))}\leq C(\|y_0\|_{L^2(\Omega)}^2+\|f\|_{L^2(Q_T)}^2+\|\Fvec\|_{\Lvec^2(Q_T)}^2).
	\end{alignedat}
\end{equation}	

	To end this second step, let us prove a uniform estimate for $\nabla y_n$. To do this, let us fix $n\geq1$ and 
	as $\nabla y_n(t)\in \Lvec^2(\Omega)$, we can write
\begin{equation} \label{def_gradientyn}
	\nabla y_n(t)= \sum_{k=1}^{+\infty}\xi_n^k(t) \uvec_k\quad \hbox{for a.e. } 0\leq t\leq T
\end{equation}
where $\xi_n^k$ denotes a time function for $[0,T]\to \mathbb{R}$ for each $k$.
Then, fixing $m \geq 1$ and multiplying the second equation of \eqref{eq:y_and_p}  by $\xi_n^k(t)$, summing over $k=1,\cdots,m$, we deduce
\[
	 (\nabla y_n(t),\sum_{k=1}^{m}\xi_n^k(t) \uvec_k) \leq{C\over2}(\|\Fvec(t)\|^2+\|\pvec_m(t)\|^2)+{1\over2}\|\sum_{k=1}^{m}\xi_n^k(t) \uvec_k\|_{\Lvec^2(\Omega)}^2.
\]
Integrating with respect to the time variable and recalling (\ref{est_p_l2}), we find 
\[
\int_0^T  (\nabla y_n(t),\sum_{k=1}^{m}\xi_n^k(t) \uvec_k) dt  \leq C\biggl(\|\Fvec(t)\|_{\Lvec^2(Q_T)}^2+\|f\|_{L^2(Q_T)}^2+\|y_0\|_{L^2(\Omega)}^2\biggr)+ {1\over2}\int_0^T\|\sum_{k=1}^{m}\xi_n^k(t) \uvec_k\|_{\Lvec^2(\Omega)}^2 dt.
\]
Let $m\to+\infty$ and using (\ref{def_gradientyn}), we finally obtain

\begin{equation}\label{est_grad_y}
	 \|\nabla y_n\|_{L^2(Q_T)}^2 \leq C(\|y_0\|_{L^2(\Omega)}^2+\|f\|_{L^2(Q_T)}^2+\|\Fvec\|_{\Lvec^2(Q_T)}^2).
\end{equation}

\item[c)] {\bf Building a weak solution.}  Let us pass the limit in the sequence $(y_n,\pvec_m)$.

	From the a priori estimates \eqref{est_y_infty}, \eqref{est_p_l2}, \eqref{est_y_t} and \eqref{est_grad_y},
	there exist subsequences $(y_{nl})_{l=1}^{\infty}\subset (y_{n})_{n=1}^{\infty}$
	and $(\pvec_{ml})_{l=1}^{\infty}\subset (\pvec_{m})_{m=1}^{\infty}$ and functions 
	$\pvec\in \Lvec^2(Q_T)$ and $y\in L^2(0,T;H_0^1(\Omega))$ with $y_t\in L^2(0,T;H^{-1}(\Omega))$ 
	such that
\begin{equation}
\label{eq:weakconvergence}
	\begin{alignedat}{2}
			&y_{nl}\to y \quad \hbox{weakly in}\quad L^2(0,T;H_0^1(\Omega)),\\
			&y_{nl}'\to y' \quad \hbox{weakly in}\quad L^2(0,T;H^{-1}(\Omega)),\\
			&\pvec_{ml}\to \pvec \quad \hbox{weakly in}\quad L^2(Q_T).
	\end{alignedat}
\end{equation}

Next, let  $w\in C^1([0,T];span\{w_k\}_{k=1}^{r})$ with $r\leq n$ and $ \uvec\in C^1([0,T];span\{\uvec_k\}_{k=1}^s)$ with
	$s\leq m$. Taking $w$ and $\uvec$ as the test functions in (\ref{eq:y_and_p}) and integrating with respect to time, we obtain
\begin{equation}\label{eq:y_and_p_intl2}
	\left\{
		\begin{array}{ll}
			\noalign{\smallskip}\dis\int_0^T[\langle y'_n,w\rangle+(\pvec_m,\nabla w) 
			+ (d\,y_n,w)]\,dt=\int_0^T(f,w)\,dt,\\
			 \noalign{\smallskip}\dis
			 \int_0^T[(\nabla y_n,\uvec)-(c^{-1}\pvec_m,\uvec)]\,dt =\int_0^T(c^{-1}\Fvec,\uvec)\,dt.
		\end{array}
		\right.
\end{equation}
 Taking $n=nl$ and $m=ml$ in the above equalities and passing to the limit, we obtain, in view of the weak convergence 
 in (\ref{eq:weakconvergence}), we also have
\begin{equation}\label{eq:y_and_p_intl2_lim}
	\left\{
		\begin{array}{ll}
			\noalign{\smallskip}\dis\int_0^T[\langle y',w\rangle+(\pvec,\nabla w) 
			+ (d\,y,w)]\,dt=\int_0^T(f,w)\,dt,\\
			 \noalign{\smallskip}\dis
			 \int_0^T[(\nabla y,\uvec)-(c^{-1}\pvec,\uvec)]\,dt =\int_0^T(c^{-1}\Fvec,\uvec)\,dt.
		\end{array}
		\right.
\end{equation}
Eventually, by a density property of $C^1([0,T];span\{w_k\}_{k=1}^{r})$ and  $C^1([0,T];span\{\uvec_k\}_{k=1}^s)$ in $L^2(0,T;H^1_0(\Omega))$ and $\Lvec^2(Q_T)$ respectively, it follows from (\ref{eq:y_and_p_intl2_lim}) that items $(i)$ and $(ii)$ of Definition \ref{def:weak_heat_stand_mixed} hold true for the pair $(y,\pvec)$.
\item[d)] {\bf Initial datum.} We now check that item $(iii)$ of Definition \ref{def:weak_heat_stand_mixed} holds true as well. First, since $y\in L^2(0,T;H_0^1(\Omega))$ with $y_t\in L^2(0,T;H^{-1}(\Omega))$
we can deduce that $y\in C^0([0,T];L^2(\Omega))$ (see Theorem $3$, pag. 303, in \cite{evans}). So, from \eqref{eq:y_and_p_intl2_lim},
	we deduce that
\begin{equation}\label{eq:y_and_p_intl2_lim2}
	\begin{array}{ll}
			\noalign{\smallskip}\dis\int_0^T[-\langle y,w'\rangle+(\pvec,\nabla w) 
			+ (d\,y,w)]\,dt=\int_0^T(f,w)\,dt+(y(\cdot,0),w(\cdot,0)),
	\end{array}
\end{equation}
	for all $w\in C^1([0,T];H^1_0(\Omega))$ such that $w(\cdot,T)=0$.
	
	And from \eqref{eq:y_and_p_intl2}, we have the same for the sequence $(y_{ml},p_{ml})$
\begin{equation}\label{eq:y_and_p_intl2_2}
	\begin{array}{ll}
			\noalign{\smallskip}\dis\int_0^T[-\langle y_{nl},w'\rangle+(\pvec_{ml},\nabla w) 
			+ (d\,y_{nl},w)]\,dt=\int_0^T(f,w)\,dt+(y_{nl}(\cdot,0),w(\cdot,0)).
	\end{array}
\end{equation}
	Taking the limit, we obtain
\begin{equation}\label{eq:y_and_p_intl2_3}
	\begin{array}{ll}
			\noalign{\smallskip}\dis\int_0^T[-\langle y,w'\rangle+(\pvec,\nabla w) 
			+ (d\,y,w)]\,dt=\int_0^T(f,w)\,dt+(y_0,w(\cdot,0)).
	\end{array}
\end{equation}

	Comparing \eqref{eq:y_and_p_intl2_lim2} and \eqref{eq:y_and_p_intl2_3}, we conclude $y(\cdot,0)=y_0$.

\item[e)] {\bf Uniqueness.} The uniqueness is deduced from the energy estimate (\ref{energy_est_mixed}).
\end{itemize}
\Fin

		It is worth mentioning that the existence and uniqueness of weak solution for the parabolic in the mixed form with $y_0\in H^1_0(\Omega)$ was proved
		in \cite{jerome_space_time} for $\Fvec=\ovec$ (Theorem $2.1$) and in \cite{li_todd_huang} for the general case
		(Lemma $3.1$).  
	%	{\color{blue}  yes, it is worth mentioning that, but in \cite{jerome_space_time,li_todd_huang}, a different $L^2-H(div)$ variational formulation is considered ! why and what is the link with the formulation of Definition A.1 ?   I mean that it should have an explanation to introduce a $L^2-H(div)$ formulation !}

	Now, let us introduce another definition of solution for \eqref{eq:heat_mixed_form}.
\begin{definition}\label{weak_heat_transposition}
	We say that the pair $(y,\pvec)\in L^2(Q_T)\times \Lvec^2(Q_T)$ is a solution by transposition of \eqref{eq:heat_mixed_form} if and only if :
 $$	
 \nonumber
 	\jjntQT (y(x,t),\pvec(x,t))\cdot( g(x,t),\Gvec(x,t))\,dx\,dt= \Mvec(g, \Gvec)\quad\forall (g,\Gvec)\in L^2(Q_T)\times \Lvec^2(Q_T)
$$
	with $\Mvec:L^2(Q_T)\times \Lvec^2(Q_T)\to \mathbb{R}$ given by
$$
\nonumber
	\Mvec(g, \Gvec):=\jjntQT f(x,t)\varphi(x,t)\,dx\,dt+(y_0,\varphi(\cdot,0))+\jjntQT c^{-1}\Fvec(x,t)\cdot\boldsymbol{\sigma}(x,t)\,dx\,dt
$$
	where $(\varphi,\boldsymbol{\sigma})$ is the unique weak solution of
\begin{equation}
\label{eq:heat_mixed_form_back}
	\left\{
		\begin{array}{lll}
   			-\varphi_t- \nabla\cdot \boldsymbol{\sigma}+d\,\varphi =g, \quad \nabla \varphi-c^{-1}\boldsymbol{\sigma}=\Gvec 	& \textrm{in}& Q_T,    		\\
   			\varphi  = 0										  		&\textrm{on}& \Sigma_T, 	\\
   			\varphi(x, T) = 0 								  		& \textrm{in}& \Om.
   		\end{array} 
 	\right.
\end{equation}
\end{definition}

\begin{remark}\label{rmq_mixed_transp}
	According to Definition \eqref{eq:heat_mixed_form}, a solution by transposition for \eqref{eq:heat_mixed_form} is a solution by transposition for the standard
	parabolic equation \eqref{eq:parab}.
	Indeed, taking $\Gvec=\ovec$, we obtain
 $$	
 \nonumber
 	\jjntQT y(x,t)g(x,t)\,dx\,dt= \jjntQT f(x,t)\varphi(x,t)\,dx\,dt+(y_0,\varphi(\cdot,0))
	+\jjntQT\Fvec(x,t)\cdot\nabla\varphi(x,t)\,dx\,dt,
$$	
	for any $g\in L^2(Q_T)$, which is the definition of solution by transposition for 
	the equation \eqref{eq:parab}.
\end{remark}
\begin{proposition}\label{ex_un_mixed_sol2}
	There exists a unique solution by transposition for \eqref{eq:heat_mixed_form}. 
\end{proposition}
\textsc{Proof-} Firstly, notice that  $ \Mvec:L^2(Q_T)\times \Lvec^2(Q_T)\to \mathbb{R}$
	is a linear form. Then, since $(\varphi,\boldsymbol{\sigma})$ is the unique weak solution, we obtain, in view of Proposition \ref{ex_un_mixed_sol2}, 
\[
	\|\varphi\|_{L^2(Q_T)}^2+\|\varphi(\cdot,0)\|_{L^2(Q_T)}^2+\|\boldsymbol{\sigma}\|_{\Lvec^2(Q_T)}^2\leq C \|(g,\Gvec)\|_{L^2(Q_T)\times \Lvec^2(Q_T)}^2.	
\]
	This implies that the linear form $\Mvec$ is continuous. Therefore, by the Riesz representation Theorem, 
	there exists a unique pair
$$
\nonumber
	(y,\pvec)\in L^2(Q_T)\times \Lvec^2(Q_T)
$$
	such that
$$	
\nonumber
 	\jjntQT (y(x,t),\pvec(x,t))\cdot( g(x,t),\Gvec(x,t))\,dx\,dt= \Mvec(g, \Gvec)\quad\forall (g,\Gvec)\in L^2(Q_T)\times \Lvec^2(Q_T).
$$

	The uniqueness is obtained by {\it Du Bois-Raymond's Lemma}.
\Fin

%{\color{blue} i think it is worth to say something about the link between the weak solution and the transposition solution.}

\end{appendix}

% \newpage   
%%%% REFERENCES

\bibliographystyle{siam}
%\bibliography{biblio_IP_heat_Diego}

%\begin{acknowledgements}
%\end{acknowledgements}


\end{document}